\documentclass[11pt,a4paper]{article}
\usepackage[utf8]{inputenc}
\usepackage[T1]{fontenc}
\usepackage{amsmath}
\usepackage{amsthm}
\usepackage{amsfonts}
\usepackage{amssymb}
\usepackage{mathtools}
\usepackage{mathrsfs}
\usepackage{url}
\usepackage{mathdots}
\usepackage[english]{babel}
\usepackage{geometry}
\usepackage{verbatim}
\usepackage{hyperref}
\usepackage{float}
\usepackage{authblk}
\usepackage{enumerate}
\usepackage{tikz-cd}
\usepackage{extarrows}
\usepackage{comment}
\usepackage{bbm}
\usepackage[style=alphabetic, backend=bibtex, backref=true]{biblatex}
\usepackage{csquotes}
\bibliography{vanEst}

\hypersetup{
    bookmarks=true,         
    unicode=false,          
    pdftoolbar=true,        
    pdfmenubar=true,        
    pdffitwindow=false,     
    pdfstartview={FitH},    
  pdftitle={A Nomizu-van Est theorem in Ekedahl's \(\ell\)-adic setting},    
    pdfauthor={Olivier Taïbi},     
    colorlinks=true,       
    linkcolor=blue,          
    citecolor=green,        
    filecolor=green,      
    urlcolor=cyan}           

\newcommand{\Q}{\mathbb{Q}}
\newcommand{\A}{\mathbb{A}}
\newcommand{\R}{\mathbb{R}}

\newcommand{\N}{\mathbb{N}}
\newcommand{\Qp}{\mathbb{Q}_p}

\newcommand{\Qell}{\mathbb{Q}_{\ell}}

\newcommand{\Zell}{\mathbb{Z}_{\ell}}
\newcommand{\Z}{\mathbb{Z}}
\newcommand{\id}{\mathrm{id}}

\newcommand{\Tot}{\operatorname{Tot}}

\newcommand{\Rep}{\operatorname{Rep}}
\newcommand{\Fil}{\operatorname{Fil}}
\newcommand{\gr}{\operatorname{gr}}
\newcommand{\Pol}{\operatorname{Pol}}
\newcommand{\pol}{\operatorname{pol}}

\newcommand{\Spec}{\operatorname{Spec}}

\newcommand{\AW}{\operatorname{AW}}
\newcommand{\EML}{\operatorname{EML}}
\mathchardef\mhyphen="2D

\newcommand{\norm}{\mathrm{norm}}

\newcommand{\Lie}{\operatorname{Lie}}

\newcommand{\im}{\operatorname{im}}

\newcommand{\mfrak}{\mathfrak{m}}
\newcommand{\nfrak}{\mathfrak{n}}

\newcommand{\alg}{\mathrm{alg}}

\newcommand{\cont}{\mathrm{cont}}
\newcommand{\car}{\operatorname{char}}
\newcommand{\fg}{\mathrm{fg}}

\newcommand{\GL}{\mathrm{GL}}

\newcommand{\End}{\mathrm{End}}
\newcommand{\Hom}{\operatorname{Hom}}

\newcommand{\res}{\operatorname{res}}
\newcommand{\ind}{\mathrm{ind}}

\newcommand{\Ocal}{\mathcal{O}}
\newcommand{\Dcal}{\mathcal{D}}

\newcommand{\Hbf}{\mathbf{H}}
\newcommand{\Gbf}{\mathbf{G}}

\newcommand{\Nbf}{\mathbf{N}}
\newcommand{\Bbf}{\mathbf{B}}

\newcommand{\Ubf}{\mathbf{U}}

\newcommand{\NvE}{\nu}

\newcommand{\ul}[1]{\underline{#1}}

\newtheorem{theo}{Theorem}[section]
\newtheorem{lemm}[theo]{Lemma}
\newtheorem{coro}[theo]{Corollary}
\newtheorem{defi}[theo]{Definition}
\newtheorem{prop}[theo]{Proposition}

\newtheorem{rema}[theo]{Remark}
\newtheorem{assu}[theo]{Assumption}

\numberwithin{equation}{section}

\begin{document}

\baselineskip=16pt

\author{Olivier Taïbi}

\title{A Nomizu-van Est theorem in Ekedahl's derived \(\ell\)-adic setting}

\maketitle

\begin{abstract}
  A theorem of Nomizu and van Est computes the cohomology of a compact nilmanifold, or equivalently the group cohomology of an arithmetic subgroup of a unipotent linear algebraic group over \(\Q\).
  We prove a similar result for the cohomology of a compact open subgroup of a unipotent linear algebraic group over \(\Qell\) with coefficients in a complex of continuous \(\ell\)-adic representations.
  We work with the triangulated categories defined by Ekedahl which play the role of ``derived categories of continuous \(\ell\)-adic representations''.
  This is motivated by Pink's formula computing the derived direct image of an \(\ell\)-adic local system on a Shimura variety in its minimal compactification, and its application to automorphic perverse sheaves on Shimura varieties.
  The key technical result is the computation of the cohomology with coefficients in a unipotent representation with torsion coefficients by an explicit complex of polynomial cochains which is of finite type.
\end{abstract}

\setcounter{tocdepth}{2}
\tableofcontents
\newpage

\section{Introduction}

A theorem of Nomizu \cite{Nomizu} and van Est \cite{vanEst} identifies, for a unipotent linear algebraic group \(\Nbf\) over \(\R\) and a discrete cocompact subgroup \(\Gamma\) of \(\Nbf(\R)\), the group cohomology groups \(H^i(\Gamma, \R)\) with the Lie algebra cohomology groups \(H^i(\Lie \Nbf(\R), \R)\).
This also holds with coefficients in a non-trivial representation of \(\Nbf\).
A rational version of this isomorphism, with \(\Nbf\) defined over \(\Q\) and \(\Gamma\) an arithmetic subgroup of \(\Nbf(\Q)\), was obtained in \cite[\S 24]{GoreskyHarderMacPherson_weighted}\footnote{More precisely, half of this isomorphism was constructed loc.\ cit., the morphism from algebraic group cohomology to group cohomology.
  For the other half, the morphism from algebraic group cohomology to Lie algebra cohomology, see Section \ref{sec:alg_gp_coh_and_Lie}.}.
For a prime number \(\ell\) this result was used\footnote{To be honest I do not understand Pink's reduction, see Remark \ref{rem:dont_understand_Pink}.} by Pink \cite[Theorem 5.3.1]{Pink_ladic_Shim} to deduce from his Theorem 4.2.1 loc.\ cit.\ the explicit computation of the restriction to strata of the higher direct images of \(\ell\)-adic automorphic local systems on a Shimura variety in its minimal compactification.
This computation also uses Kostant's formula \cite[Theorem 5.14]{Kostant_Liealgcoh} which explicitly computes Lie algebra cohomology in this setting.
Pink's Theorem 4.2.1 identifies the restriction to strata of higher direct images of local systems with torsion coefficients with local systems obtained from certain group cohomology groups.
In fact this theorem is stated in derived categories, and for the study of certain perverse sheaves on minimal compactifications of Shimura varieties it is desirable to have a ``derived'' analogue of \cite[Theorem 5.3.1]{Pink_ladic_Shim}, that is a commutative square between suitable triangulated categories.
Such an analogue seems to have been implicitly used in \cite[Théorème 2.2.1]{MorelSiegel1}, \cite[Theorem 1.2.3]{MorelBook} and \cite{Zhu_IHorth}.
To prove such an analogue one needs the \(\ell\)-adic analogue of \cite[Theorem 4.2.1]{Pink_ladic_Shim}, which we have not considered seriously yet but should be rather formal to deduce from the proof loc.\ cit., and an \(\ell\)-adic analogue of the Nomizu-van Est theorem in suitable triangulated categories (``derived categories of continuous \(\ell\)-adic representations of profinite groups'').
The goal of this article is to prove such an analogue, Theorem \ref{thm:ladic_vanEst}.

The rough strategy is familiar: for a unipotent linear algebraic group \(\Nbf\) over \(\Qell\), an algebraic representation \(V\) of \(\Nbf\) and a compact open subgroup \(N\) of \(\Nbf(\Qell)\) we would like to define two natural quasi-isomorphisms
\begin{equation} \label{eq:intro_strategy}
  C^\bullet_{\cont}(N, V) \longleftarrow C^\bullet_{\pol}(\Nbf, V) \longrightarrow C^\bullet_{\mathrm{CE}}(\Lie \Nbf, V)
\end{equation}
where \(C^\bullet_{\cont}\) denotes the complex of continuous cochains, \(C^\bullet_{\pol}\) denotes the complex of polynomial cochains, the left map is the obvious one, \(C^\bullet_{\mathrm{CE}}\) denotes the Chevalley-Eilenberg complex computing Lie algebra cohomology, and the right map is obtained by taking differentials at the origin (details in Section \ref{sec:alg_gp_coh_and_Lie}).
One could also replace polynomial cochains by analytic cochains, in fact this gives quasi-isomorphisms in much greater generality  \cite[Theorems V.2.3.10 and V.2.4.10]{Lazard_gpanp}.
Unfortunately in the triangulated categories that will consider this strategy will not work as is, mainly because the objects that we need to consider are not these complexes.
Indeed continuous group cohomology with \(\ell\)-adic coefficients does not quite fit in the derived functor formalism, and instead we work with Ekedahl's formalism introduced in \cite{Ekedahl_adic}.
This formalism applies to arbitrary toposes and seems well-suited for the above application which mixes continuous representations of profinite topological groups and sheaves on schemes.
In particular this formalism defines, for a profinite topological group \(K\) and a finite extension \(E\) of \(\Qell\), a certain triangulated category \(D^+(K, E)\), which admits a natural functor \(F\) from the derived bounded category of finite-dimensional continuous representations of \(K\) over \(E\) (see Corollary \ref{cor:functor_F}).
This seems to be the only way to explicitly construct objects in \(D^+(K,E)\).
We emphasize that \(D^+(K,E)\) is only defined after the integral analogue \(D^+(K, \Ocal_E)\) where \(\Ocal_E\) is the ring of integers of \(E\), essentially by inverting \(\ell\), and that objects of \(D^+(K, \Ocal_E)\) are inverse systems of complexes of smooth \(\Ocal_E/\mfrak_E^i[K]\)-modules where \(\mfrak_E\) is the maximal ideal of \(\Ocal_E\).
In particular the case of torsion coefficients seems to be difficult to avoid when proving results about such objects.
For \(N\) a closed normal subgroup of \(K\) we have by general topos theory \cite[\S 5]{Ekedahl_adic} a derived pushforward functor
\[ R\Gamma(N, -): D^+(K, E) \longrightarrow D^+(K/N, E) \]
and our goal is to compute the composition \(R\Gamma(N, F(-))\) in certain cases, again in terms of bounded complexes of finite-dimensional representations of \(K/N\) over \(E\).
The problem is that the middle term in \eqref{eq:intro_strategy}, i.e.\ the complex of polynomial cochains, does not make sense a priori as an object of \(D^+(K/N, E)\), because its terms have infinite dimension over \(E\).
Even if we try to make sense of it, it still does not seem possible to define the right map in \eqref{eq:intro_strategy} over \(\Ocal_E\), even after multiplying by a large power of \(\ell\).
Our solution roughly consists of replacing the complex of polynomial cochains by an explicit subcomplex which is quasi-isomorphic to it and whose terms have finite dimension over \(E\).
For the purpose of induction and to keep an action of \(K/N\) it is crucial that these subcomplexes be canonical, even if they depend on auxiliary suitable integers, and that we obtain results for representations with torsion coefficients first.

Let us describe the contents of the paper.
As above let \(E\) be a finite extension of \(\Qell\), with ring of integers \(\Ocal_E\) having maximal ideal \(\mfrak_E\).
In Sections \ref{sec:Ekedahl} and \ref{sec:gp_coh_derived} we review the setting in which we will work, Ekedahl's formalism for (discrete) sets with continuous action of a profinite group, and gather a few lemmas.
In section \ref{sec:gp_coh_cont_cochains} we show that \(R\Gamma(N, -)\) can be computed using the usual complex of locally constant cochains in towers in \(\Ocal_E/\mfrak_E^i\)-modules, at least in the case where \(K\) is a semi-direct product \(N \rtimes H\).
For the purpose of induction we recall in section \ref{sec:expl_HS} the explicit ``Künneth formalism'' in the simplicial setting, due to Eilenberg, Zilber and MacLane, for the universal cover of the classifying space of a semi-direct product of groups.
The resulting Künneth formula is useless because the associated topological spaces are contractible, but by equivariance this formalism yields an explicit version of the Hochschild-Serre spectral sequence in the case of semi-direct products which is robust when imposing additional conditions on cochains.
To start the induction we consider in Section \ref{sec:pol_cochains_Zell} the case where \(N \simeq \Zell\), which is enlightened by the explicit description of \(1\)-cocycles and motivates the introduction of filtrations on \(\Ocal_E[N]\)-modules.
We show in Lemma \ref{lem:Cpol_Zell_cplx} that we can compute \(R\Gamma(\Zell, -)\) on bounded below complexes of unipotent \(\Ocal_E/\mfrak_E^i[\Zell]\)-modules using polynomial cochains of bounded degree (in terms of a filtration).
In Section \ref{sec:pol_cochains_general} we consider the general case where \(N\) is a torsion-free nilpotent topologically finitely generated pro-\(\ell\) group (equivalently, \(N\) is a compact open subgroup of \(\Nbf(\Qell)\) for some unipotent linear algebraic group \(\Nbf\) over \(\Qell\), see Corollary \ref{cor:all_N_sat_assu}).
It is natural to use Mal'cev coordinates to describe polynomial functions on \(N\) but to obtain a canonical notion of degree we are led to introduce ponderations, which unfortunately complicate the filtrations.
We complete the induction for torsion coefficients in Theorem \ref{thm:main_torsion}, pass to Ekedahl's triangulated categories in Corollary \ref{cor:derived_integral_main_result}, and generalize this result suitably after extension of scalars from \(\Ocal_E\) to \(E\) in Proposition \ref{pro:cont_coh_as_pol_ell_inv}.
This last result is where we get rid of the assumption that \(K\) is a semi-direct product \(N \rtimes H\).
We recall, in the form suitable for our application, the definitions of Lie algebra cohomology and Hochschild cohomology for algebraic groups in Sections \ref{sec:Lie_alg_coh} and \ref{sec:alg_gp_coh_and_Lie}.
In Section \ref{sec:pf_vanEst} we put the pieces together and finally prove Theorem \ref{thm:ladic_vanEst}.

I am grateful to Sophie Morel for her comments on a previous version of this paper.

\section{``Derived categories'' of \(\ell\)-adic representations}
\label{sec:Ekedahl}

Let \(E\) be a finite extension of \(\Qell\), \(\Ocal_E\) its ring of integers, \(\mfrak_E\) its maximal ideal.
Let \(K\) be a topological group.
We consider the category \(S_K\) of sets\footnote{in a fixed universe} with a smooth action of \(K\), i.e.\ the stabilizer of any point is an open.
It is a topos: Giraud's axioms \cite[Exposé IV Théorème 1.2]{SGA4-1} are readily verified (colimits and finite limits are computed as in the category of sets, with obvious action of \(K\)).

\begin{rema}
  To avoid giving a false sense of generality, we point out that denoting by \(K^0\) the connected component of \(1 \in K\), every smooth action of \(K\) factors through the totally disconnected topological group \(K/K^0\), and so \(S_K\) is isomorphic to \(S_{K/K^0}\).
  So we might as well assume that \(K\) is totally disconnected.
  Also recall that any locally compact totally disconnected topological group is locally profinite (van Dantzig's theorem).

  In fact in this paper we are mostly concerned with the case where \(K\) is profinite (in which case the topos \(S_K\) was considered in \cite[Exposé IV \S 2.7.2]{SGA4-1}).
\end{rema}

For \(? \in \{+,b\}\) denote by \(D^?(K, \Ocal_E)\) the associated ``derived category of \(\Ocal_E\)-modules with continuous action of \(K\)'' \cite{Ekedahl_adic}.
We borrow notation from loc.\ cit., in particular we will need to work with the ringed topos \((S_K^{\N}, (\Ocal_E)_\bullet)\).
Note that condition B loc.\ cit.\ holds.

\begin{prop} \label{pro:F_to_Eke}
  Let \(\Rep_{\fg, \cont}(K, \Ocal_E)\) be the abelian category of finitely generated \(\Ocal_E\)-modules with a continuous linear action of \(K\).
  The functor
  \begin{align*}
    F: \Rep_{\fg, \cont}(K, \Ocal_E) & \longrightarrow (S_K^{\N}, (\Ocal_E)_\bullet)-\text{modules} \\
    A & \longmapsto (A/\mfrak_E^i A)_{i \geq 1}
  \end{align*}
  induces an exact functor \(F: D^b(\Rep_{\fg, \cont}(K, \Ocal_E)) \to D^b(K, \Ocal_E)\).
\end{prop}
\begin{proof}
  The preservation of cones is obvious.
  By Lemma 3.2 loc.\ cit.\ the functor \(F\) maps complexes concentrated in one degree to \(\Ocal_E\)-complexes.
  By an immediate induction it maps bounded complexes to objects of \(D^b(K, \Ocal_E)\).

  We are left to check that \(F\) maps exact complexes to negligible complexes.
  By Lemma 1.4 loc.\ cit.\ it is enough to check that the image of an exact complex is essentially zero, which easily follows from the Artin-Rees lemma for \(\Ocal_E\)-modules.
\end{proof}

We will often consider images \(F(A^\bullet)\) of this functor \(F\) where \(A^\bullet\) is a bounded below complex of objects of \(\Rep_{\fg, \cont}(K,\Ocal_E)\) for which there exists \(n_0 \in \Z\) such that we have \(H^n(A^\bullet) = 0\) for any \(n > n_0\).
We emphasize that in this case \(F(A^\bullet)\) is defined as the bounded complex \(F(\tau_{\leq n_0} A^\bullet)\), which up to a well-defined isomorphism in \(D^b(K, \Ocal_E)\) does not depend on the choice of \(n_0\).

\begin{lemm}
  Let \((\Dcal, \Sigma)\) be an \(\Ocal_E\)-linear triangulated category\footnote{The class of distinguished triangles is abusively omitted from the notation.}.
  Let \(\Dcal[\ell^{-1}]\) be the category having the same objects as \(\Dcal\), with
  \[ \Hom_{\Dcal[\ell^{-1}]}(X, Y) := E \otimes_{\Ocal_E} \Hom_{\Dcal}(X, Y). \]
  The auto-equivalence \(\Sigma\) of \(\Dcal\) induces an auto-equivalence of \(\Dcal[\ell^{-1}]\).
  Call a triangle in \(\Dcal[\ell^{-1}]\) distinguished if it is isomorphic to the image, under the obvious functor \(\Dcal \to \Dcal[\ell^{-1}]\), of a distinguished triangle of \(\Dcal\).
  Then \(\Dcal[\ell^{-1}]\) is a triangulated category.
\end{lemm}
\begin{proof}
  We check the axioms of \cite[Ch.\ II Définition 1.1.1]{Verdier_catder}.
  We first check axiom TR1.
  The first part follows immediately from the definition.
  The fact that \(X \to X \to 0 \to \Sigma X\) is distinguished for any object \(X\) is obvious.
  Any morphism \(X \to Y\) in \(\Dcal[\ell^{-1}]\) is of the form \(\ell^{-k} f\) for some integer \(k \geq 0\) and some morphism \(f : X \to Y\) in \(\Dcal\).
  We have a commutative diagram in \(\Dcal[\ell^{-1}]\)
  \[ \begin{tikzcd}
    X \arrow[r, "{f}"] \arrow[d, "{1}"] & Y \arrow[d, "{\ell^{-k}}"] \\
    X \arrow[r, "{\ell^{-k} f}"] & Y
  \end{tikzcd} \]
  where the vertical morphism are isomorphisms, and the existence of a distinguished triangle containing \(\ell^{-k} f\) in \(\Dcal[\ell^{-1}]\) follows from the existence of a distinguished triangle containing \(f\) in \(\Dcal\).
  The last part of TR1 is tautological.

  Axiom TR2 is formal.

  To check axiom TR3 we can reduce to the case of a diagram
  \[ \begin{tikzcd}
    X \arrow[r, "{u}"] \arrow[d, "{\ell^{-a} f}"] & Y \arrow[d, "{\ell^{-b} g}"]
    \arrow[r, "{v}"] & Z \arrow[r, "{w}"] & \Sigma X \arrow[d, "{\ell^{-a}
    \Sigma f}"] \\
    X' \arrow[r, "{u'}"] & Y' \arrow[r, "{v'}"] & Z' \arrow[r, "{w'}"] & \Sigma
    X'
  \end{tikzcd} \]
  where \(u,v,w,u',v',w',f,g\) are morphisms in \(\Dcal\) and both rows are distinguished triangles in \(\Dcal\).
  Letting \(k=\max(a,b)\) and using the commutative diagram
  \[ \begin{tikzcd}
    X \arrow[r, "{u}"] \arrow[d, "{\ell^k}"] & Y \arrow[d, "{\ell^k}"] \arrow[r,
    "{v}"] & Z \arrow[d, "{\ell^k}"] \arrow[r, "{w}"] & \Sigma X \arrow[d,
    "{\ell^k}"] \\
    X \arrow[r, "{u}"] & Y \arrow[r, "{v}"] & Z \arrow[r, "{w}"] & \Sigma X
  \end{tikzcd} \]
  where all vertical morphisms are isomorphisms in \(\Dcal[\ell^{-1}]\), we reduce to the case where \(a=b=0\), and apply TR3 for \(\Dcal\).

  In order to check the octahedral axiom TR4 we may reduce, thanks to TR3, to the case where the three initial distinguished triangles are obtained (in a strict sense) from distinguished triangles in \(\Dcal\).
\end{proof}

\begin{rema}
  One can also show that \(\Dcal[\ell^{-1}]\) is equivalent to the quotient of \(\Dcal\) by the full subcategory consisting of all objects \(X\) of \(\Dcal\) for which there exists an integer \(k \geq 0\) satisfying \(\ell^k = 0\) in \(\Hom_{\Dcal}(X,X)\), and apply \cite[Theorem 10.2.3]{KashiwaraShapira_catandsheaves}.
\end{rema}

\begin{lemm} \label{lem:triang_cat_inv_ell}
  Let \(D^b(\Rep_{\fg, \cont}(K, E))\) be the abelian category of finite-dimensional vector spaces over \(E\) endowed with a continuous action of \(K\) \emph{which admit a \(K\)-stable \(\Ocal_E\)-lattice} (this condition is automatically satisfied if \(K\) is compact).
  The functor
  \begin{align*}
    D^b(\Rep_{\fg, \cont}(K, \Ocal_E))[\ell^{-1}] & \longrightarrow
    D^b(\Rep_{\fg, \cont}(K, E)) \\
    A & \longmapsto E \otimes_{\Ocal_E} A
  \end{align*}
  is an equivalence of triangulated categories.
\end{lemm}
\begin{proof}
  The fact that this functor is exact is obvious.
  Essential surjectivity follows from the fact that for a finite free \(\Ocal_E\)-module \(\Lambda\) any compact subset of \(E \otimes_{\Ocal_E} \Lambda\) is contained in some \(\mfrak_E^{-n} \Lambda\).

  Let us show that the functor is conservative.
  By exactness it is enough to show that an object \(A\) of \(D^b(\Rep_{\fg, \cont}(K, \Ocal_E))[\ell^{-1}]\) is zero if its image is zero.
  This last condition means that the cohomology groups of \(A\) are killed by \(\ell^k\) for some integer \(k \geq 0\).
  There exist integers \(a \leq b\) such that \(H^i(A) = 0\) for \(i \not\in [a,b]\), and we proceed by induction on \(b-a\).
  If \(b-a=0\) then \(A\) is quasi-isomorphic to a complex concentrated in one degree, and its endomorphism \(\ell^k\) is both invertible and zero, so \(A\) is a zero object.
  If \(b-a>0\) we use the distinguished triangle in \(D^b(\Rep_{\fg,\cont}(K,\Ocal_E))\)
  \[ H^a(A)[-a] \to A \to B \xrightarrow{+1} \]
  where \(H^i(B)=0\) for \(i \not\in [a+1,b]\).
  By induction hypothesis both \(H^a(A)\) and \(B\) are zero, so \(A\) is zero too.

  Let \(A,B\) be objects of \(D^b(\Rep_{\fg,\cont}(K, \Ocal_E))\).
  It is clear that every morphism of complexes \(E \otimes_{\Ocal_E} A \to E \otimes_{\Ocal_E} B\) is of the form \(\ell^{-k} f\) for some morphism of complexes \(f: A \to B\).
  This fact and conservativity imply fullness.

  To check faithfulness it is enough to prove that if \(f: A \to B\) is a morphism between bounded complexes of objects in \(\Rep_{\fg,\cont}(K, \Ocal_E)\) mapping to the zero morphism in \(D^b(\Rep_{\fg,\cont}(K, E))\) then it already maps to the zero morphism in \(D^b(\Rep_{\fg,\cont}(K, \Ocal_E))[\ell^{-1}]\).
  The condition means that there exists a quasi-isomorphism \(g_E: E \otimes_{\Ocal_E} B \to C_E\) of bounded complexes of objects in \(\Rep_{\fg,\cont}(K, E)\) such that the composition \(E \otimes_{\Ocal_E} A \to C_E\) is homotopic to zero.
  There exists a bounded complex \(C\) of objects in \(\Rep_{\fg,\cont}(K, \Ocal_E)\), an isomorphism \(E \otimes_{\Ocal_E} C \simeq C_E\) and a morphism of complexes \(g: B \to C\) inducing \(g_E\) such that the composition \(g \circ f: A \to C\) is homotopic to zero.
  We have proved above that \(g\) becomes an isomorphism in \(D^b(\Rep_{\fg,\cont}(K, \Ocal_E))[\ell^{-1}]\), and so \(f\) becomes the zero morphism in this category.
\end{proof}

Denote \(D^+(K, E) = D^+(K, \Ocal_E)[\ell^{-1}]\), and similarly for \(D^b\).

\begin{coro} \label{cor:functor_F}
  The functor \(F\) of Proposition \ref{pro:F_to_Eke} induces a functor, abusively still denoted by \(F\):
  \[ F: D^b(\Rep_{\fg,\cont}(K, E)) \longrightarrow D^b(K, E). \]
\end{coro}

\section{Group cohomology as a derived functor}
\label{sec:gp_coh_derived}

For any morphism of totally disconnected topological groups \(f: G \to H\) we have a morphism of toposes \(S_G \to S_H\) for which pullback is the restriction functor \(f^*\) (more commonly denoted by \(\res_G\) or \(\res^H_G\)) and the pushforward functor \(f_*\) is smooth induction \(\ind_G^H\).
More precisely for an object \(X\) of \(S_G\) we define \(\ind_G^H X\) as the set of functions \(\varphi: H \to X\) such that for any \(g \in G\) and any \(h \in H\) we have \(\varphi(f(g) h) = g \cdot \varphi(h)\) and which are invariant under some open subgroup of \(H\).
Note that in the case where \(H\) is the quotient of \(G\) by a closed normal subgroup \(N\) the pushforward functor is really the functor of ``\(N\)-invariants'', simply denoted \(-^N\).
In this case we denote by \(R\Gamma(N, -)\) its derived functor \(D^+(S_G^{\N}, (\Ocal_E)_\bullet) \to D^+(S_{G/N}^{\N}, (\Ocal_E)_\bullet)\) (resp.\ \(D^+(G, \Ocal_E) \to D^+(G/N, \Ocal_E)\), resp.\ \(D^+(G, E) \to D^+(G/N, E)\); see \cite[\S 5]{Ekedahl_adic}).
Our ultimate goal is to compute, in some cases, the composition \(R\Gamma(N, F(-))\) as the composition of \(F\) with a more concrete functor.

\begin{lemm} \label{lem:res_preserves_inj}
  Let \(K'\) be a closed subgroup of a topological group \(K\).
  Assume that \(K'\) is open in \(K\) or that \(K\) is profinite.
  Then the restriction functor
  \[ \res_{K'}: (S_K^{\N}, (\Ocal_E)_\bullet)-\text{modules} \longrightarrow (S_{K'}^{\N}, (\Ocal_E)_\bullet)-\text{modules} \]
  maps injective objects to injective objects.
\end{lemm}
\begin{proof}
  If \(K'\) is an open subgroup of \(K\) then \(f^* = \res_{K'}\) admits a \emph{left} adjoint, namely the functor \(\mathrm{find}_{K'}^K\) of induction with finite support in \(K' \backslash K\).
  More precisely the unit for this adjunction is given by the obvious identification of an \((\Ocal_E/\mfrak_E^m)_\bullet\)-module \(A_m\) in \(S_{K'}\) with the subobject of \(\res_{K'} \mathrm{find}_{K'}^K A_m\) consisting of all functions supported on \(K'\).
  The functor \(\mathrm{find}_{K'}^K\) is obviously exact so the lemma follows from this adjunction in this case.

  Now consider the case where \(K\) is profinite.
  Let \(I_\bullet\) be an injective \((\Ocal_E)_\bullet\)-module in \(S_K^{\N}\).
  Consider an embedding \(A_\bullet \subset B_\bullet\) of \((\Ocal_E)_\bullet\)-modules in \(S_{K'}^{\N}\) and a morphism \(\varphi: A_\bullet \to \res_{K'} I_\bullet\).
  We wish to extend \(\varphi\) to \(B_\bullet\).
  Assume that \(B_\bullet\) is not equal to \(A_\bullet\), i.e.\ there exists \(n \geq 1\) such that we have \(A_n \subsetneq B_n\).
  Choose \(x \in B_n \smallsetminus A_n\).
  For \(m \leq n\) let \(C_m\) be the sub-\(\Ocal_E/\mfrak_E^m[K']\)-module of \(B_m\) generated by the image of \(x\), and for \(m>n\) let \(C_m=0\).
  This defines a sub-\((\Ocal_E)_\bullet\)-module of \(B_\bullet\), and to extend \(\varphi\) to \(A_\bullet+C_\bullet\) it is enough to extend \(\varphi|_{A_\bullet \cap C_\bullet}\) to \(C_\bullet\).
  There exists an open and normal subgroup \(U\) of \(K\) such that \(x\) is fixed by \(U \cap K'\).
  In particular each \(C_m\) is an \(\Ocal_E/\mfrak_E^m[K'/(U \cap K')]\)-module.
  Since \(K'/(U \cap K')\) is finite and \(\Ocal_E\) is noetherian each \(A_m \cap C_m\) is a finitely generated \(\Ocal_E/\mfrak_E^m[K'/(U \cap K')]\)-module.
  This property (and \(C_m=0\) for \(m>n\)) imply that up to replacing \(U\) by a smaller normal open subgroup of \(K\) we have \(\varphi(A_m \cap C_m) \subset I_m^U\) for all \(m\).
  Using the isomorphism \(K'/(U \cap K') \simeq K'U/U\) we may see \(C_\bullet\) as a \((S_{K'U}^\N, (\Ocal_E)_\bullet)\)-module, \(A_\bullet \cap C_\bullet\) as a submodule and \(\varphi|_{A_\bullet \cap C_\bullet}\) as a morphism to \(\res_{K'U} I_\bullet\).
  Now \(K'U\) is an open subgroup of \(K\) so by the previous case we can extend \(\varphi|_{A_\bullet \cap C_\bullet}\) to \(C_\bullet\).
  So \(\varphi\) can be extended to give a morphism \(A_\bullet + C_\bullet \to \res_{K'} I_\bullet\).
  By the usual argument using Zorn's lemma we conclude that any \(\varphi\) as above can be extended to give a morphism \(B_\bullet \to \res_{K'} I_\bullet\).
\end{proof}

Assume for the rest of this section that the assumption of Lemma \ref{lem:res_preserves_inj} is satisfied.
Let \(N\) be a closed normal subgroup of \(K\).
Denote \(N' = K' \cap N\), a closed normal subgroup of \(K'\).
We have a morphism of composite functors
\[ \begin{tikzcd}[column sep=5em]
  D^+(S_K^{\N}, (\Ocal_E)_\bullet) \arrow[r, "{R\Gamma(N, -)}" above] \arrow[d,
  "{\res_{K'}}" left] & D^+(S_{K/N}^{\N}, (\Ocal_E)_\bullet) \arrow[d,
  "{\res_{K'/N'}}" right] \arrow[dl, Rightarrow, shorten=8mm] \\
  D^+(S_{K'}^{\N}, (\Ocal_E)_\bullet) \arrow[r, "{R\Gamma(N', -)}" below] &
  D^+(S_{K'/N'}^{\N}, (\Ocal_E)_\bullet),
\end{tikzcd} \]
i.e.\ a morphism of functors
\begin{equation} \label{eq:res_RGamma}
  r_{K,N,K'}: \res_{K'/N'}( R\Gamma(N, -) ) \longrightarrow R\Gamma(N', \res_{K'}(-)).
\end{equation}
Indeed, before derivation we have an obvious morphism of functors
\[ \res_{K'/N'}(-^N) \longrightarrow (\res_{K'} -)^{N'}. \]
By general formalism, the right derived functor of the left-hand side is isomorphic to \(\res_{K'/N'}(R\Gamma(N, -))\) because \(\res_{K'/N'}\) is exact, and the right derived functor of the right-hand side is isomorphic to \(R\Gamma(N', \res_{K'}(-))\) thanks to Lemma \ref{lem:res_preserves_inj}.

If \(K'\) contains \(N\), i.e.\ if \(N'\) is equal to \(N\), then \(r_{K,N,K'}\) is an isomorphism of functors, essentially because it is an isomorphism before derivation.

If \(K''\) is a closed subgroup of \(K'\), denoting \(N'' = K'' \cap N\) we have an equality of morphisms of functors
\begin{equation} \label{eq:comp_res_transf}
  r_{K,N,K''} = r_{K',N',K''}(\res_{K'}) \circ \res_{K''/N''}(r_{K,N,K'}).
\end{equation}
The tedious but unsurprising verification of this formula is left to the reader.

\section{Group cohomology using continuous cochains}
\label{sec:gp_coh_cont_cochains}

Now assume that \(K\) is a semi-direct product \(N \rtimes H\) of profinite topological groups.
For \(A^\bullet\) a complex of \((S_K^{\N}, (\Ocal_E)_\bullet)\)-modules consider the double complex \(C^\bullet(N, A^\bullet)\) of \((S_H^{\N}, (\Ocal_E)_\bullet)\)-modules defined by
\[ C^n(N, A^m)_i := C^n(N, A^m_i) := \{ \text{locally constant } \varphi: N^{n+1} \to A^m_i \}^N \]
where the action of \(N\) is defined by the formula
\begin{equation} \label{eq:action_N_cochains}
  (g \cdot \varphi)(x_0, \dots, x_n) = g \cdot \varphi(g^{-1} x_0, \dots, g^{-1} x_n)
\end{equation}
and with obvious morphisms \(C^n(N, A^m)_{i+1} \to C^n(N, A^m)_i\) and \(C^n(N, A^m)_i \to C^n(N, A^{m+1})_i\).
(Of course we define \(C^n(N, -) = 0\) for \(n<0\)).
The action of \(H\) on \(C^n(N, A_i^m)\) is defined by the formula
\begin{equation} \label{eq:action_H_cochains}
  (h \cdot \varphi)(x_0, \dots, x_n) = h \cdot \varphi(h^{-1} x_0 h, \dots, h^{-1} x_n h).
\end{equation}
The differentials are defined by the formula
\begin{equation} \label{eq:diff_hom_cochains}
  d \varphi(x_0, \dots, x_{n+1}) = \sum_{j=0}^{n+1} (-1)^j \varphi(x_0, \dots, \widehat{x_j}, \dots, x_{n+1})
\end{equation}
where the hat means that \(x_j\) is omitted.
It will be useful to keep in mind that these definitions come from the realization of the classifying space \(BN\) as a quotient \(N \backslash EN\) where \(EN\) is the simplicial set with \((EN)_n = N^{n+1}\), degeneracy maps for \(0 \leq i \leq n-1\)
\begin{align*}
  s_i: N^n & \longrightarrow N^{n+1} \\
  (x_0, \dots, x_{n-1}) & \longmapsto (x_0, \dots, x_i, x_i, \dots, x_{n-1}),
\end{align*}
face maps for \(0 \leq i \leq n\)
\begin{align*}
  \partial_i: N^{n+1} & \longrightarrow N^n \\
  (x_0, \dots, x_n) & \longmapsto (x_0, \dots, \widehat{x_i}, \dots, x_n)
\end{align*}
and action of \(N\) on \(EN\)
\[ g \cdot (x_0, \dots, x_n) := (g x_0, \dots, g x_n). \]
In particular it will be useful to view functions from \(N^{n+1}\) to an abelian group \(M\) as morphisms of abelian groups \(\Z[N^{n+1}] \to M\).
Letting \(\partial: \Z[N^{n+2}] \to \Z[N^{n+1}]\) be \(\sum_{i=0}^{n+1} (-1)^i \partial_i\) we have \(d \varphi = \varphi \circ \partial\).
The proof of the Dold-Kan correspondence applied to the simplicial abelian group \(\Z[EN]\) yields the decompositions
\[ \Z[N^{n+1}] = \bigcap_{0 \leq i <n} \ker \partial_i \oplus \sum_{0 \leq i < n} \im s_i \]
which are compatible with \(\partial\), and the corresponding projections \(\Z[N^{n+1}] \to \bigcap_{0 \leq i <n} \ker \partial_i\) are explicitly given by \((1 - s_{n-1} \partial_{n-1}) \dots (1 - s_0 \partial_0)\).
Defining \(\theta: \Z[N^{n+1}] \to \Z[N^{n+2}]\) for \(n \geq 0\) as
\begin{equation} \label{eq:def_theta}
  \sum_{i=0}^n (-1)^i s_i (1 - s_{i-1} \partial_{i-1}) \dots (1 - s_0 \partial_0)
\end{equation}
(and \(\theta: 0 \to \Z[N]\) in the only possible way) we have a chain homotopy relation
\[ \id_{\Z[N^{n+1}]} = (1 - s_{n-1} \partial_{n-1}) \dots (1 - s_0 \partial_0) + \partial \theta + \theta \partial. \]
If \(M\) is an \(\Ocal_E[K]\)-module we have a corresponding decomposition of complexes of \(\Ocal_E[H]\)-modules
\[ C^\bullet(N, M) = C^\bullet_{\norm}(N, M) \oplus C^\bullet_{\deg}(N, M) \]
where
\begin{align}
  C^n_{\norm}(N, M) &:= \left\{ \varphi \in C^n(N, M) \,\middle|\, \forall 0 \leq i < n, \, \varphi \circ s_i = 0 \right\} \label{eq:def_norm_cochains} \\
  C^n_{\deg}(N, M) &:= \left\{ \varphi \in C^n(N, M) \,\middle|\, \varphi|_{\cap_{0 \leq i < n} \ker \partial_i} = 0 \right\}. \nonumber
\end{align}
and denoting
\begin{align} \label{eq:norm_cochains_homotopy}
  \theta^* : C^{n+1}(N, M) & \longrightarrow C^n(N, M) \\
  \varphi & \longmapsto \varphi \circ \theta \nonumber
\end{align}
the corresponding projection \(C^n(N, M) \to C^n_{\deg}(N, M)\) is equal to \(d \circ \theta^* + \theta^* \circ d\).
In particular the inclusion of complexes \(C^\bullet_{\norm}(N, M) \to C^\bullet(N, M)\) is a quasi-isomorphism.
This is of course well-known but we will need these results for certain sub-complexes of cochains defined by additional conditions on cochains which will be obviously preserved by all the maps above.

\begin{lemm} \label{lem:derived_gp_coh_as_cocycles}
  Let \(K\) be a semi-direct product \(N \rtimes H\) of profinite topological groups.
  \begin{enumerate}
  \item If \(A^\bullet \to B^\bullet\) is a quasi-isomorphism between bounded-below complexes of \((S_K^{\N}, (\Ocal_E)_\bullet)\)-modules then the induced morphism \(\Tot^\bullet(C^\bullet(N, A^\bullet)) \to \Tot^\bullet(C^\bullet(N, B^\bullet))\) is a quasi-isomorphism.
    In other words we have a well-defined functor
    \begin{align*}
      R\Gamma_{\cont}(N, -): D^+(S_K^{\N}, (\Ocal_E)_\bullet) & \longrightarrow D^+(S_{K/N}^{\N}, (\Ocal_E)_\bullet) \\
      A^\bullet & \longmapsto \Tot^\bullet(C^\bullet(N, A^\bullet)).
    \end{align*}
    This functor is exact.
  \item The functors \(R\Gamma(N, -)\) and \(R\Gamma_{\cont}(N, -)\) are isomorphic.
  \end{enumerate}
\end{lemm}
\begin{proof}
  The functors \(C^n(N, -)\) are exact so the preservation of quasi-isomorphisms in the first point follows from \cite[Theorem 12.5.4]{KashiwaraShapira_catandsheaves}.
  The fact that the functor \(R\Gamma_{\cont}(N, -)\) is exact, i.e.\ maps cones to cones, is a routine computation with double complexes.

  The morphism of functors \(R\Gamma(N,-) \to R\Gamma_{\cont}(N, -)\) is defined as follows.
  For an object \(A^\bullet\) of \(D^+(S_K^{\N}, (\Ocal_E)_\bullet)\) choose a quasi-isomorphism \(A^\bullet \to I^\bullet\) where \(I^n = 0\) for \(n<<0\) and each \(I^n\) is an injective \((S_K^{\N}, (\Ocal_E)_\bullet)\)-module.
  This induces a quasi-isomorphism \(R\Gamma_{\cont}(N, A^\bullet) \to R\Gamma_{\cont}(N, I^\bullet)\).
  We have a natural morphism of double complexes \((I^\bullet)^N \to C^\bullet(N, I^\bullet)\) where \((I^\bullet)^N\) abusively denotes the double complex whose \(0\)-th row is \((I^\bullet)^N\) and whose other rows are zero.
  By definition we have \(R\Gamma(N, A^\bullet) = (I^\bullet)^N\).
  The composition in \(D^+(S_H^{\N}, (\Ocal_E)_\bullet)\)
  \[ R\Gamma(N, A^\bullet) = (I^\bullet)^N \longrightarrow R\Gamma_{\cont}(N, I^\bullet) \xleftarrow{\sim} R\Gamma_{\cont}(N, A^\bullet) \]
  is functorial in \(A^\bullet\), as usual because \(D^+(S_K^{\N}, (\Ocal_E)_\bullet)\) is equivalent to the homotopy category of bounded below chain complexes consisting of injective \((S_K^{\N}, (\Ocal_E)_\bullet)\)-modules.

  To check that this morphism of functors is an isomorphism it is enough to check that for any injective \((S_K^{\N}, (\Ocal_E)_\bullet)\)-module \(I\) the complex
  \[ I^N \longrightarrow C^\bullet(N, I) \]
  is a resolution of the \((S_H^{\N}, (\Ocal_E)_\bullet)\)-module \(I^N\).
  Thanks to Lemma \ref{lem:res_preserves_inj} we may assume \(K=N\).
  For any normal open subgroup \(U\) of \(N\) the \((S_{N/U}^{\N}, (\Ocal_E)_\bullet)\)-module \(I^U\) is injective.
  For any integer \(n \geq 0\) and any integer \(i \geq 1\) we have an identification
  \[ C^n(N/U, I_i^U) = \Hom_{(S_{N/U}^{\N}, (\Ocal_E)_\bullet)}(R(i,n), I^U) \]
  where \(R(i,n)\) is the \((S_{N/U}, (\Ocal_E)_\bullet)\)-module defined by
  \[ R(i,n)_j = \begin{cases}
    \Ocal_E/\mfrak_E^j[(N/U)^{n+1}] & \text{ if } 1 \leq j \leq i \\
    0 & \text{ if } j>i
  \end{cases} \]
  with obvious transition morphisms.
  Defining an \((S_{N/U}, (\Ocal_E)_\bullet)\)-module \(R(i)\) similarly:
  \[ R(i)_j = \begin{cases}
    \Ocal_E/\mfrak_E^j & \text{ if } 1 \leq j \leq i \\
    0 & \text{ if } j>i,
  \end{cases} \]
  it is well-known that the complex
  \[ \dots \to R(i, 2) \to R(i, 1) \to R(i, 0) \to R(i) \to 0 \]
  is exact, and we obtain, by almost the same argument as in the case of cohomology for finite groups, that \(C^\bullet(N/U, I^U)\) is a resolution of \((I^U)^{N/U} = I^N\).
  We conclude by taking the direct limit over all normal open subgroups \(U\) of \(N\).
\end{proof}

\begin{lemm} \label{lem:cmp_res_inj_cocyc}
  Let \(K\) be a profinite topological group.
  Let \(N\) be a closed normal subgroup of \(K\).
  Let \(K'\) be a closed subgroup of \(K\) and denote \(N' = K' \cap N\).
  The following diagram of morphisms of functors is commutative.
  \[ \begin{tikzcd}
    \res_1(R\Gamma(N, -)) \arrow[r, "{\res_1(r_{K,N,K'})}"] \arrow[d, "{r_{K,N,N}}" right, "{\sim}" left] & \res_1(R\Gamma(N', \res_{K'}(-))) \arrow[d, "{r_{K',N',N'}}" right, "{\sim}" left] \\
    R\Gamma(N, \res_N(-)) \arrow[r, "{r_{N,N,N'}}"] \arrow[d, "{\sim}"] & R\Gamma(N', \res_{N'}(-)) \arrow[d, "{\sim}"] \\
    R\Gamma_{\cont}(N, \res_N(-)) \arrow[r, "{\res}"] & R\Gamma_{\cont}(N', \res_{N'}(-))
  \end{tikzcd} \]
  For the top horizontal morphism, the identification \(\res_1 \simeq \res_1 \circ \res_{K'/N'}\) was used implicitly.
  The bottom vertical isomorphisms are the ones defined in Lemma \ref{lem:derived_gp_coh_as_cocycles}.
  The bottom horizontal morphism is the one induced by the usual ``restriction of cocycles'' morphism of functors \(C^\bullet(N, -) \to C^\bullet(N', \res_{N'}(-))\).
\end{lemm}
\begin{proof}
  Commutativity of the top square is a special case of \eqref{eq:comp_res_transf}, so we are left to prove commutativity of the bottom square.
  
  We may restrict to considering bounded-below complexes \(I^\bullet\) consisting of injective \((S_N^{\N}, (\Ocal_E)_\bullet)\)-modules, in which case we have \(R\Gamma(N, I^\bullet) = (I^\bullet)^N\) and also
  \[ R\Gamma(N', \res_{N'}(I^\bullet)) = (\res_{N'}(I^\bullet))^{N'} \]
  thanks to Lemma \ref{lem:res_preserves_inj}.
  The following diagram of double complexes, where as in the previous proof complexes are considered as double complexes concentrated in one row, is clearly commutative.
  \[ \begin{tikzcd}
    (I^\bullet)^N \arrow[d] \arrow[r] & (\res_{N'} I^\bullet)^{N'} \arrow[d] \\
    C^\bullet_{\cont}(N, I^\bullet) \arrow[r] & C^\bullet_{\cont}(N', \res_{N'}(I^\bullet))
  \end{tikzcd} \]
  Taking total complexes yields the result.
\end{proof}

\section{Explicit Hochschild-Serre for semi-direct products}
\label{sec:expl_HS}

By standard composition of derived pushforward functors between toposes, if \(N_1 \subset N_2\) are closed normal subgroups of a locally profinite group \(K\) then the morphism of functors \(R\Gamma(N_2, -) \to R\Gamma(N_2/N_1, R\Gamma(N_1, -))\) is an isomorphism.
We will need an explicit, more flexible version of this statement for \(R\Gamma_{\cont}\).
Assume that \(K\) decomposes as a semi-direct product \(N_1 \rtimes (N_2 \rtimes H)\).
The comparison between \(R\Gamma_{\cont}(N_2, R\Gamma_{\cont}(N_1, -))\) and \(R\Gamma_{\cont}(N_1 \rtimes N_2, -)\) proceeds from the proof of the Künneth formula for simplicial sets which is due to Eilenberg-MacLane.
More precisely we have an obvious identification \(E(N_1 \rtimes N_2) \simeq EN_1 \times EN_2\) and the resulting action of \(N_1 \rtimes N_2\) on the right-hand side is simple: \(N_1\) acts as usual on the factor \(EN_1\) and trivially on the factor \(EN_2\) while \(N_2\) acts ``by conjugation'' on the factor \(EN_1\) and as usual on the factor \(EN_2\).
Of course the Künneth formula for \(EN_1 \times EN_2\) is not very interesting because both factors are contractible, but we will see that since all maps are very explicit it is almost trivial to add the ``invariance under \(N_1 \rtimes N_2\)'' condition everywhere.

For \(i,j \in \Z_{\geq 0}\) let \(\Sigma(i,j)\) denote the set of \((j,i)\)-shuffles, i.e.\ the set of pairs \((\sigma, \tau)\) of non-decreasing maps \(\sigma: [0,i+j] \to [0,j]\) and \(\tau: [0,i+j] \to [0,i]\) satisfying \(\sigma(a+1)+\tau(a+1)-\sigma(a)-\tau(a)=1\) for all \(0 \leq a < i+j\).
Note that this condition implies that both \(\sigma\) and \(\tau\) are surjective, and that \(\Sigma(i,j)\) is naturally in bijection with the set of subsets of \([0,i+j[\) of cardinality \(j\) (resp.\ \(i\)), via
\[ (\sigma, \tau) \longmapsto I_\sigma := \left\{ 0 \leq a < i+j \,\middle|\, \sigma(a+1)>\sigma(a) \right\} \]
(resp.\ idem with \(\tau\) replacing \(\sigma\)).
For \((\sigma,\tau) \in \Sigma(i,j)\) define its sign \(\epsilon(\sigma, \tau)\) as \((-1)^x\) where \(x\) is the cardinality of the set
\[ \left\{ (a,b) \in I_\sigma \times I_\tau \,\middle|\, a>b \right\}. \]
For \(M\) a smooth \(\Ocal_E[N_1 \rtimes N_2]\)-module define the Alexander-Whitney and Eilenberg-MacLane maps and a map \(\Psi\) as follows:
\begin{align*}
  \AW&: \bigoplus_{i+j=n} C^i(N_2, C^j(N_1, M)) \longrightarrow C^n(N_1 \rtimes N_2, M) \\
  \EML&: C^n(N_1 \rtimes N_2, M) \longrightarrow \bigoplus_{i+j=n} C^i(N_2, C^j(N_1, M)) \\
  \Psi&: C^{n+1}(N_1 \rtimes N_2, M) \longrightarrow C^n(N_1 \rtimes N_2, M)
\end{align*}
by the formulas
\begin{equation} \label{eq:def_AW}
  \AW((\varphi_{i,j})_{i,j})(x_0 y_0, \dots, x_n y_n) = \sum_{i+j=n} \varphi_{i,j}(y_j, \dots, y_n)(x_0, \dots, x_j),
\end{equation}
\begin{equation} \label{eq:def_EML}
  \EML(\varphi)_{i,j}(y_0, \dots, y_i)(x_0, \dots, x_j) = \sum_{(\sigma, \tau) \in \Sigma(i,j)} \epsilon(\sigma, \tau) \varphi((x_{\sigma(a)} y_{\tau(a)})_{0 \leq a \leq i+j}),
\end{equation}
\begin{multline} \label{eq:def_Psi}
  \Psi(\varphi)(x_0 y_0, \dots, x_n y_n) = \\
  \sum_{p+q<n} \sum_{(\sigma, \tau) \in \Sigma(q,p+1)} (-1)^m \epsilon(\sigma, \tau) \varphi((x_a y_a)_{0 \leq a \leq m}, (x_{m+\sigma(a)} y_{n-q+\tau(a)})_{0 \leq a \leq p+q+1})
\end{multline}
where in the sum \(m=n-p-q-1\).
It is easy to check that these maps are well-defined, i.e.\ they preserve the smoothness and invariance conditions on cochains.

\begin{theo} \label{thm:Eilenberg_Zilber}
  \begin{enumerate}
  \item The operator \(\AW\) maps
    \[ \bigoplus_{i+j=n} C^i_{\norm}(N_2, C^j_{\norm}(N_1, M)) \text{ to } C^n_{\norm}(N_1 \rtimes N_2, M), \]
    and vice-versa for \(\EML\).
  \item The maps \(\AW\) and \(\EML\) commute with differentials, where on the total complex \((\bigoplus_{i+j=n} C^i(N_2, C^j(N_1, M)))_n\) the differentials are defined by summing the maps
    \[ d_{N_1}^j: C^i(N_2, C^j(N_1, M)) \to C^i(N_2, C^{j+1}(N_1, M)) \]
    and
    \[ (-1)^j d_{N_2}^i: C^i(N_2, C^j(N_1, M)) \to C^{i+1}(N_2, C^j(N_1, M)). \]
  \item The restriction of \(\EML \circ \AW\) to \(\bigoplus_{i+j=n} C^i_{\norm}(N_2, C^j_{\norm}(N_1, M))\) is the identity.
    In particular there is an explicit (see previous section) homotopy between \(\EML \circ \AW\) and \(\id\).
  \item
    On \(C^n(N_1 \rtimes N_2, M)\) we have
    \[ \AW \circ \EML - \id = d \circ \Psi + \Psi \circ d. \]
  \end{enumerate}
  In particular \(\AW\) and \(\EML\) are both quasi-isomorphisms, inverse of each other after taking cohomology.
\end{theo}
\begin{proof}
  The first point is easy and left to the reader.
  This other points are dual to the Eilenberg-Zilber theorem \cite[Theorem 2.1a]{EilenbergMacLane2}, following \cite[\S 5]{EilenbergMacLane1}.
  Our operator \(\Psi\) is dual to the one defined on p.55 of \cite{EilenbergMacLane2}.
  The closed formula for \(\Psi\) is a slight reformulation of the dual to \emph{the opposite of} the ``Shih operator'' \cite[Définition 3.1.1]{Rubio_these}.
  It is not difficult to check that the opposite of Rubio's \(\mathrm{SHI}\) satisfies the formula defining Eilenberg-MacLane's \(\Psi\) by induction.
  Details are left to the reader.
\end{proof}

We will in fact not need the precise formula for \(\Psi\) in the sequel, only its ``shape'', from which it will be obvious that it preserves certain additional conditions on cochains.

\section{Polynomial cochains for \(\Zell\)}
\label{sec:pol_cochains_Zell}

Let \(A\) be a finitely generated \(\Ocal_E/\mfrak_E^i\)-module (for some integer \(i \geq 1\)), endowed with a continuous action of \(\Zell\).
Let \(A\) be a torsion abelian group endowed with a continuous action of \(\Zell\).
We have a canonical decomposition \(A = A_\ell \times A'\) where each element of \(A_\ell\) is killed by some power of \(\ell\) and each element of \(A'\) is killed by some positive integer not divisible by \(\ell\).
Denote by \(\gamma\) the image of \(1 \in \Zell\) in \(\operatorname{Aut}(A)\).
It is well-known (see \cite[VIII.4 and XIII.1]{Serre_corpsloc}) that we have a quasi-isomorphism
\begin{equation} \label{eq:cont_coh_Zell_general}
  [A_\ell \xrightarrow{\gamma-1} A_\ell ] \longrightarrow{} C^\bullet(\Zell, A)
\end{equation}
where the complex on the left is concentrated in degrees \(0\) and \(1\) and the maps are
\begin{align*}
  A_\ell & \longrightarrow C^0(\Zell, A) \\
  a & \longmapsto (x \mapsto \gamma^x \cdot a)
\end{align*}
and
\begin{align*}
  A_\ell & \longrightarrow Z^1(\Zell, A) \\
  a & \longmapsto \left( (x,y) \mapsto \sum_{k \geq 1} \binom{y-x}{k} (\gamma-1)^{k-1} \gamma^x \cdot a \right)
\end{align*}
where the sum on the right only has finitely many non-zero terms because the endomorphism \(\gamma-1\) of \(A_\ell\) is locally nilpotent.
This formula gives the unique continuous (homogeneous) \(1\)-cocycle whose value at \((0,1)\) is \(a\).
Motivated by this simple case, we will for certain \(\Ocal_E/\mfrak_E^i\)-modules with continuous action of \(\Zell\) (and more generally, complexes) find a bounded quasi-isomorphic subcomplex of the complex of continuous cochains which consists of finitely generated \(\Ocal_E\)-modules.

For an \(\Ocal_E\)-module \(A\) and an integer \(r \geq 0\) let \(\Pol(\Zell^r, A)\) be the \(\Ocal_E\)-module of functions \(\Zell^r \to A\) of the form
\begin{equation} \label{eq:pol_func}
  (x_0, \dots, x_{r-1}) \mapsto \sum_{\ul{i} = (i_0, \dots, i_{r-1}) \in \Z_{\geq 0}^r} \binom{x_0}{i_0} \dots \binom{x_{r-1}}{i_{r-1}} a_{\ul{i}}
\end{equation}
where \(a_{\ul{i}} \in A\) vanish for almost all \(\ul{i}\).
By a well-known argument the \(a_{\ul{i}}\)'s can be recovered successively by evaluating at tuples of non-negative integers.
We say that such a function has total degree \(\leq d\) if \(a_{\ul{i}} = 0\) whenever \(|\ul{i}| > d\), where \(|\ul{i}| := i_0 + \dots + i_{r-1}\).

\begin{lemm}
  Let \(N\) be a topological group isomorphic to \(\Zell^r\) for some \(r \in \Z_{\geq 0}\).
  The notions of polynomial function \(N \to A\) and of total degree do not depend on the choice of isomorphism \(N \simeq \Zell^r\).
\end{lemm}
\begin{proof}
  Two choices of isomorphisms \(N \simeq \Zell^r\) differ by an element of \(\GL_r(\Zell)\).
  The lemma follows from Vandermonde's identity and the expansion
  \begin{equation} \label{eq:change_var_pol_Hilb}
    \binom{\lambda x}{k} = \lambda^k \binom{x}{k}+ \sum_{0 \leq i \leq k-1} a_i(\lambda, k) \binom{x}{i}
  \end{equation}
  with \(a_i(\lambda, k) \in \Zell\), valid for any \(\lambda \in \Zell\) and \(k \in \Z_{\geq 0}\).
  Details are left to the reader.
\end{proof}

Let \(N\) be a topological group isomorphic to \(\Zell^r\) for some \(r \geq 0\) and let \(A\) be an \(\Ocal_E[N]\)-module.
Let \((\Fil^j A)_{j \in \Z}\) be a filtration on \(A\) by sub-\(\Ocal_E[N]\)-modules satisfying \(\Fil^j A \subseteq \Fil^{j+1} A\) for all \(j \in \Z\) and such that \(\Fil^j A = 0\) for some \(j \in \Z\).
Our convention in this paper is that ``filtration'' always means ``non-decreasing filtration''.
For the purpose of induction we have to introduce an integer \(w>0\).
We will frequently assume that \(N\) acts trivially on each ``\(w\)-graded piece'' \(\Fil^j A / \Fil^{j-w} A\).
For any integer \(n \geq 0\) we have an action of \(N\) on \(\Pol(N^{n+1}, A)\) defined by the formula (with multiplicative notation)
\[ (g \cdot f)(h_0, \dots, h_n) = g \cdot f(g^{-1} h_0, \dots, g^{-1} h_n). \]
We endow \(\Pol(N^{n+1}, A)\) with the filtration for which \(\Fil^j_w \Pol(N^{n+1}, A)\) is the sub-\(\Ocal_E\)-module of \(\Pol(N^{n+1}, A)\) consisting of all polynomials of the form
\[ (x_0, \dots, x_n) \mapsto \sum_{i \in I} P_i(x_0, \dots, x_n) a_i \]
where \(I\) is a finite set and for which there exists a family \((d_i, j_i)_{i \in I}\) such that \(P_i\) has total degree \(\leq d_i\) and we have \(a_i \in \Fil^{j_i} A\) and \(w d_i + j_i \leq j\).
In other words, choosing an isomorphism \(N \simeq \Zell^r\) we see that \(f \in \Pol(N^{n+1}, A)\) given by coefficients \(a_{\ul{i}}\) as in \eqref{eq:pol_func} belongs to \(\Fil^j_w \Pol(N^{n+1}, A)\) if and only if for all \(\ul{i} \in \Z_{\geq 0}^{r(n+1)}\) we have \(a_{\ul{i}} \in \Fil^{j-|\ul{i}|} A\).
For \(n,m \geq 0\) we have an obvious identification of filtered \(\Ocal_E[N]\)-modules \(\Pol(N^{n+1}, \Pol(N^{m+1}, A)) \simeq \Pol(N^{n+m+2},A)\).

\begin{defi}
  Let \((A, \Fil^\bullet A)\) be a filtered \(\Ocal_E[N]\)-module satisfying \(\Fil^j A = 0\) for some \(j \in \Z\).
  For an integer \(n \geq 0\) let \(C^n_{\pol}(N, A)\) be the sub-\(\Ocal_E\)-module of \(N\)-invariant elements in \(\Pol(N^{n+1}, A)\), endowed with the usual differentials and the filtration \(\Fil^\bullet_w\) induced by that on \(\Pol(N^{n+1}, A)\).
\end{defi}

\begin{lemm} \label{lem:Fil_Cpol_Fil_Pol_Zell}
  Let \((A, \Fil^\bullet A)\) be a filtered \(\Ocal_E[\Zell]\)-module satisfying \(\Fil^j A = 0\) for some \(j \in \Z\).
  Assume that each \(\Fil^j A\) is a finitely generated \(\Ocal_E\)-module with continuous action of \(\Zell\) for the \(\mfrak_E\)-adic topology, and that \(\Zell\) acts trivially on each \(\Fil^j A / \Fil^{j-w} A\).
  Then for any \(d \in \Z\) we have an isomorphism of \(\Ocal_E\)-modules
  \begin{align*}
    \Fil^d_w C^n_{\pol}(\Zell, A) & \longrightarrow \Fil^d_w \Pol(\Zell^n, A) \\
    \varphi & \longmapsto \left( (x_1, \dots, x_n) \mapsto \varphi(0, x_1, \dots, x_n) \right).
  \end{align*}
\end{lemm}
\begin{proof}
  To avoid confusion we will denote the action of \(\Zell\) on \(A\) as
  \begin{align*}
    \Zell \times A & \longrightarrow A \\
    (x, a) & \longmapsto \gamma^x \cdot a.
  \end{align*}
  The fact that the morphism in the lemma is well-defined is clear, and injectivity immediately follows from the equivariance condition which reads
  \[ \varphi(x_0, \dots, x_n) = \gamma^{x_0} \cdot \varphi(0, x_1-x_0, \dots, x_n-x_0). \]
  We are left to check surjectivity.
  Let \(\psi \in \Fil^d_w \Pol(\Zell^n, A)\).
  The polynomial function
  \[ (x_0, \dots, x_n) \mapsto \psi(x_1-x_0, \dots, x_n-x_0) \]
  clearly belongs to \(\Fil^d_w \Pol(\Zell^{n+1}, A)\), and so there is an expansion
  \[ \psi(x_1-x_0, \dots, x_n-x_0) = \sum_{i \in I} P_i(x_0, \dots, x_n) a_i \]
  where \(I\) is a finite set, \(P_i \in \Pol(\Zell^{n+1}, \Ocal_E)\) has total degree \(\leq d_i\) and \(a_i \in \Fil^{d-wd_i} A\).
  Fix \(j_0 \in \Z\) for which we have \(\Fil^{j_0} A = 0\).
  For \(x_0 \in \Z_{\geq 0}\) we have
  \begin{align*}
    \gamma^{x_0} \cdot \psi(x_1-x_0, \dots, x_n-x_0)
    &= \sum_{i \in I} P_i(x_0, \dots, x_n) \gamma^{x_0} \cdot a_i \\
    &= \sum_{i \in I} \sum_{0 \leq k < (d-j_0)/w-d_i} P_i(x_0, \dots, x_n) \binom{x_0}{k} (\gamma-1)^k a_i
  \end{align*}
  where the last equality follows from the usual binomial identity for
  \[ \gamma^{x_0} = ((\gamma-1) + 1)^{x_0} \in \End(A) \]
  and the fact that \((\gamma-1)^k\) maps \(\Fil^{d-wd_i} A\) to \(\Fil^{d-wd_i-wk} A\).
  By continuity of the action this equality holds for arbitrary \(x_0 \in \Zell\) and it is clear on this expression that
  \[ (x_0, \dots, x_n) \longmapsto \gamma^{x_0} \cdot \psi(x_1-x_0, \dots, x_n-x_0) \]
  belongs to \(\Fil^d_w \Pol(\Zell^{n+1}, A)\).
\end{proof}

\begin{lemm} \label{lem:seed_pol_cocy}
  Let \((A,\Fil^\bullet A)\) be a filtered \(\Ocal_E[\Zell]\)-module.
  Assume that there exists \(j \in \Z\) for which we have \(\Fil^j A = 0\), that for any \(j \in \Z\) the \(\Ocal_E\)-module \(\Fil^jA\) is finitely generated and that the action of \(\Zell\) on it is continuous, and that the action of \(\Zell\) on each \(\Fil^j A / \Fil^{j-w} A\) is trivial.
  For any \(d \in \Z\) we have a quasi-isomorphism
  \[ [\Fil^d A \xrightarrow{\gamma-1} \Fil^{d-w} A] \to \Fil^d_w C^\bullet_{\pol}(\Zell, A) \]
  where the complex on the left is concentrated in degrees \(0\) and \(1\).
  In particular if \(A\) is killed by \(\mfrak_E^i\) for some integer \(i \geq 1\) and if we have \(\Fil^{d-w} A = A\) then the embedding of complexes
  \[ \Fil^d_w C^\bullet_{\pol}(\Zell, A) \to C^\bullet(\Zell, A) \]
  is a quasi-isomorphism.
\end{lemm}
\begin{proof}
  To avoid confusion we will denote the action of \(\Zell\) on \(A\) as
  \begin{align*}
    \Zell \times A & \longrightarrow A \\
    (x, a) & \longmapsto \gamma^x \cdot a.
  \end{align*}
  Let \(\sigma_w \Pol(\Zell, A)\) be \(\Pol(\Zell, A)\) but with shifted filtration:
  \[ \Fil^j \sigma_w \Pol(\Zell, A) := \Fil^{j-w}_w \Pol(\Zell, A). \]
  We use the resolution (short exact sequence)
  \[ 0 \to A \to \Pol(\Zell, A) \to \sigma_w \Pol(\Zell, A) \to 0 \]
  where the first map is \(a \mapsto (x \mapsto a)\) and the second map is \(\Delta: f \mapsto (x \mapsto f(x+1) - f(x))\).
  Note that both morphisms are strict, in fact we have a section \(s\) of \(\Delta\):
  \[ \left( x \mapsto \binom{x}{i} a \right) \mapsto \left( x \mapsto \binom{x}{i+1} a \right) \]
  which identifies \(\Pol(\Zell, A)\) with \(A \oplus \sigma_w \Pol(\Zell, A)\) as a filtered \(\Ocal_E\)-module, but \(s\) is not \(\Zell\)-equivariant.
  We get for each integer \(n \geq 0\) a short exact sequence
  \[ 0 \to \Pol(\Zell^{n+1}, A) \to \Pol(\Zell^{n+1}, \Pol(\Zell, A)) \to \Pol(\Zell^{n+1}, \sigma_w \Pol(\Zell, A)) \to 0 \]
  in which both morphisms are strict.
  By the previous lemma taking \(\Fil^d\) and \(\Zell\)-invariants yields a short exact sequence
  \begin{equation} \label{eq:short_ex_seq_alg_cocycles}
    0 \to C^n_{\pol, \leq d}(\Zell, A) \to C^n_{\pol, \leq d}(\Zell, \Pol(\Zell, A)) \to C^n_{\pol, \leq d}(\Zell, \sigma_w \Pol(\Zell, A)) \to 0
  \end{equation}

  Next we claim that the complex
  \[ C^\bullet_{\pol, \leq d}(\Zell, \Pol(\Zell, A)) \]
  is exact in positive degree.
  Using the identification of filtered \(\Ocal_E[N]\)-modules
  \begin{align*}
    \Pol(\Zell^{n+1}, \Pol(\Zell, A)) & \longrightarrow \Pol(\Zell^{n+2}, A) \\
    \varphi & \longmapsto \left( (x_0, \dots, x_{n+1}) \mapsto \varphi(x_1, \dots, x_{n+1})(x_0) \right)
  \end{align*}
  we obtain
  \[ \Fil^d_w C^n_{\pol}(\Zell, \Pol(\Zell, A)) \simeq \Fil^d_w C^{n+1}_{\pol}(\Zell, A) \]
  and applying the previous lemma to the right hand side gives us an isomorphism
  \begin{align*}
    \Fil^d_w C^n_{\pol}(\Zell, \Pol(\Zell, A)) & \longrightarrow \Fil^d_w \Pol(\Zell^{n+1}, A) \\
    \varphi & \longmapsto \left( (x_0, \dots, x_n) \mapsto \varphi(x_0, \dots, x_n)(0) \right).
  \end{align*}
  Via this identification the exactness of \(C^\bullet_{\pol, \leq d}(\Zell, \Pol(\Zell, A))\) in positive degree is obtained as usual using the homotopy corresponding to
  \begin{align*}
    \Fil^d \Pol(\Zell^{n+2}, A) & \longrightarrow \Fil^d \Pol(\Zell^{n+1}, A) \\
    \psi & \longrightarrow \left( (x_0, \dots, x_n) \mapsto \psi(0, x_0, \dots, x_n) \right).
  \end{align*}
  In degree \(0\) we find that the kernel of \(\Fil^d_w C^0_{\pol}(\Zell, \Pol(\Zell, A)) \to \Fil^d_w C^1_{\pol}(\Zell, \Pol(\Zell, A))\) is naturally isomorphic to \(\Fil^d A\).

  The same arguments are valid with \(\Pol(\Zell, A)\) replaced by \(\sigma_w \Pol(\Zell, A)\) (e.g.\ by shifting the filtration on \(A\)).
  Using the short exact sequences \eqref{eq:short_ex_seq_alg_cocycles} and a familiar double complex argument we obtain that the complexes \(\Fil^d_w C^\bullet_{\pol}(\Zell, A)\) and
  \begin{equation} \label{eq:FildA_gm1_FildmwA}
    \Fil^d A \xlongrightarrow{\gamma-1} \Fil^{d-w} A.
  \end{equation}
  are both quasi-isomorphic to the same complex.
  By making this explicit or by direct computation now that we know that \(\Fil^d_w C^\bullet_{\pol}(\Zell, A)\) is exact in degree \(>1\) we find that the complex \(\tau_{\leq 1} \Fil^d_w C^\bullet_{\pol}(\Zell, A)\) is isomorphic to \eqref{eq:FildA_gm1_FildmwA} via the isomorphisms
  \begin{align*}
    \Fil^d_w C^0_{\pol}(\Zell, A) & \longrightarrow \Fil^d A \\
    \varphi & \longmapsto \varphi(0)
  \end{align*}
  and
  \begin{align*}
    \Fil^d_w Z^1_{\pol}(\Zell, A) & \longrightarrow \Fil^{d-w} A \\
    \varphi & \longmapsto \varphi(0, 1).
  \end{align*}
  The last assertion of the lemma follows using the quasi-isomorphism \eqref{eq:cont_coh_Zell_general} recalled above.
\end{proof}

For the purpose of induction the following direct consequence of the previous lemma will be more useful.

\begin{lemm} \label{lem:Cpol_Zell_cplx}
  Let \((A^\bullet, \Fil^\bullet A^\bullet)\) be a bounded below filtered complex of \(\Ocal_E[\Zell]\)-modules.
  Assume that there exists \(j \in \Z\) such that the complex \(\Fil^j A^\bullet\) is zero, that for any \(j,n \in \Z\) the \(\Ocal_E\)-module \(\Fil^j A^n\) is finitely generated and the action of \(\Zell\) on it is continuous, and that the action of \(\Zell\) on each \(\Fil^j A^n / \Fil^{j-w} A^n\) is trivial.
  \begin{enumerate}
  \item Assume that each \(\Fil^j A^n\) is torsion, i.e.\ killed by \(\mfrak_E^i\) for some integer \(i \geq 1\), and that there exists \(d_0 \in \Z\) such that for any \(j \geq d_0\) the embedding of complexes \(\Fil^j A^\bullet \to \Fil^{j+1} A^\bullet\) is a quasi-isomorphism.
    Then for any integers \(d \geq d_0+w\) and \(j \geq d\) the embedding of complexes
    \[ \Tot^\bullet(\Fil^d_w C_{\pol}^\bullet(\Zell, A^\bullet)) \to \Tot^\bullet(C^\bullet(\Zell, \Fil^j A^\bullet)) \]
    is a quasi-isomorphism.
  \item Let \(d \in \Z\) be such that \(\Fil^d A^\bullet\) and \(\Fil^{d-w} A^\bullet\) are exact in degree \(>r\).
    Then \(\Tot^\bullet(\Fil^d_w C_{\pol}^\bullet(\Zell, A^\bullet))\) is exact in degree \(>r+1\).
  \end{enumerate}
\end{lemm}
\begin{proof}
  \begin{enumerate}
  \item By Lemma \ref{lem:seed_pol_cocy} the embedding of double complexes
    \begin{equation} \label{eq:dbl_cplxes_FCpolZell}
      \left[ \Fil^d A^\bullet \xrightarrow{\gamma-1} \Fil^{d-w} A^\bullet \right] \to \Fil^d_w C_{\pol}^\bullet(\Zell, A^\bullet)
    \end{equation}
    where the left hand side has two rows in degrees \(0\) and \(1\), consists of quasi-isomorphisms on columns.
    Because we assume \(j \geq d \geq d_0+w\) the embedding of double complexes
    \[ \left[ \Fil^d A^\bullet \xrightarrow{\gamma-1} \Fil^{d-w} A^\bullet \right] \to \left[ \Fil^j A^\bullet \xrightarrow{\gamma-1} \Fil^j A^\bullet \right] \]
    consists of quasi-isomorphisms on all rows.
    We also know that the embedding of double complexes
    \[ \left[ \Fil^j A^\bullet \xrightarrow{\gamma-1} \Fil^j A^\bullet \right] \to C^\bullet(\Zell, \Fil^j A^\bullet) \]
    consists of quasi-isomorphisms on all columns (this is where we use the assumption that each \(\Fil^j A^n\) is torsion).
    These three embeddings and the embedding
    \[ \Fil^d_w C^\bullet_{\pol}(\Zell, A^\bullet) \to C^\bullet(\Zell, \Fil^j A^\bullet) \]
    sit in an obvious commutative diagram, and the result follows.
  \item This also follows from the fact that the embedding of double complexes \eqref{eq:dbl_cplxes_FCpolZell} consists of quasi-isomorphisms on columns.
\end{enumerate}
\end{proof}

\section{Polynomial cochains for torsion-free pro-\(\ell\) nilpotent groups}
\label{sec:pol_cochains_general}

\begin{lemm} \label{lem:char_tors_free_nilp}
  Let \(N\) be a pro-\(\ell\) topological group, i.e.\ \(N\) is isomorphic to a projective limit of finite \(\ell\)-groups.
  The following conditions are equivalent.
  \begin{enumerate}
  \item \label{it:tors_free_nilp}
    \(N\) is nilpotent, torsion-free and topologically finitely generated.
  \item \label{it:central_series_Zell_ri}
    There exists a filtration
    \[ N=N_0 \supset N_1 \supset \dots \supset N_m \supset N_{m+1} = 0 \]
    where each \(N_i\) is a closed normal subgroup in \(N\), each \(N_i/N_{i+1}\) is isomorphic to \(\Zell^{r_i}\) for some \(r_i \in \Z_{\geq 0}\), and the action by conjugation of \(N\) on \(N_i/N_{i+1}\) is trivial.
  \end{enumerate}
\end{lemm}
\begin{proof}
  The fact that \eqref{it:central_series_Zell_ri} implies \eqref{it:tors_free_nilp} is easy and left to the reader.
  For the other implication we claim that we can take \((N_i)_i\) to be the upper central series of \(N\).
  Let \((N_j')_j\) be the lower central series of \(N\).
  Applying \cite[Lemma 17.2.1]{KargapolovMerzljakov_GTM} to finite quotients of \(N\) we see that each \(N'_j\) is also topologically finitely generated, and so we have isomorphisms \(N'_j/N'_{j+1} \simeq \Zell^{r_j'} \times \Z/\ell^{s'_j}\Z\) for all \(j\).
  Note that any closed subgroup of \(\Zell^{r_j'} \times \Z/\ell^{s'_j}\Z\) is also topologically finitely generated.
  Intersecting \((N'_j)_j\) with \(N_i\) then shows that each \(N_i\) is also topologically finitely generated, and so we have isomorphisms \(N_i/N_{i+1} \simeq \Zell^{r_i} \times \Z/\ell^{s_i}\Z\).
  We are left to show that all \(s_i\) vanish, i.e.\ that \(N_i/N_{i+1}\) is \(\ell\)-torsion-free.
  This can be proved by descending induction on \(i \in \{0, \dots, m\}\), using the well-known fact that for any \(i \in \{0, \dots, m-1\}\), \(x \in N_i\) and \(y \in N\) we have \([x^\ell,y] = [x,y]^\ell\) in \(N_{i+1}/N_{i+2}\).
\end{proof}

Let \((N,(N_i)_i)\) be as in Lemma \ref{lem:char_tors_free_nilp} \eqref{it:central_series_Zell_ri}.
If we choose \(\ul{g} = (g_{i,j})_{0 \leq i \leq m, 1 \leq j \leq r_i}\) where \(g_{i,j} \in N_i\) and the image of \((g_{i,j})_{1 \leq j \leq r_i}\) is a basis of the \(\Zell\)-module \(N_i/N_{i+1}\), we obtain a homeomorphism between \(\Zell^{r_m + \dots + r_0}\) and \(N\),
\begin{equation} \label{eq:def_psi_g}
  \psi_{\ul{g}}: (x_{m,1}, \dots, x_{m,r_m}, \dots, x_{0,1}, \dots, x_{0,r_0}) \longmapsto g_{m,1}^{x_{m,1}} \dots g_{0,r_0}^{x_{0,r_0}}.
\end{equation}
The \(x_{i,j}\)'s may be called Mal'cev coordinates after \cite{Malcev_nilp_tors_free} where the case of discrete torsion-free nilpotent groups was considered.
A simple argument by induction shows that for any other choice \(\ul{h}\) of lifts of bases of \((N_i/N_{i+1})_i\) the composition \(\psi_{\ul{h}}^{-1} \circ \psi_{\ul{g}}\) is of the form
\begin{equation} \label{eq:change_basis_fil_N}
  (\ul{x}_m, \dots, \ul{x}_0) \longmapsto (\alpha_m(\ul{x}_m) + \beta_m(\ul{x}_{m-1}, \dots, \ul{x}_0), \dots, \alpha_1(\ul{x}_1) + \beta_1(\ul{x}_0),\alpha_0(\ul{x}_0))
\end{equation}
where for each \(i\) we have \(\alpha_i \in \GL_{r_i}(\Zell)\) and \(\beta_i: \Zell^{r_{i-1} + \dots + r_0} \to \Zell^{r_i}\) is a continuous map (which depends on \(\ul{g}\) and \(\ul{h}\))
Similarly for any family \(\ul{g}\) as above and any \(h \in N\) we have
\[ h^{-1} \psi_{\ul{g}}(\ul{x}_m, \dots, \ul{x}_0) = \psi_{\ul{g}}(\ul{x}_m + P_m(\ul{x}_{m-1}, \dots, \ul{x}_0), \dots, \ul{x}_0 + P_0) \]
where each \(P_i: \Zell^{r_{i-1} + \dots + r_0} \to \Zell^{r_i}\) is a continuous map (which depends on \(\ul{g}\) and \(h\)).

To obtain an intrinsic definition of polynomials with some notion of bounded degree, we will require the following assumption, but see Lemma \ref{lem:alg_unip_satisfies_assu} and Remark \ref{rem:assu_always_sat} below.

\begin{assu} \label{assu:N_alg}
  There exists an integer \(d\) (which depends on \(N\) and \((N_i)_i\)) such that the following holds.
  \begin{enumerate}
  \item For any families \(\ul{g}\) and \(\ul{h}\) as above, each \(\beta_i\) is a polynomial of degree \(\leq d\).
  \item For any family \(\ul{g}\) as above and any \(h \in N\), each \(P_i\) is a polynomial of degree \(\leq d\).
  \end{enumerate}
\end{assu}

\begin{lemm} \label{lem:alg_unip_satisfies_assu}
  Let \(\Nbf\) be a unipotent algebraic group over \(\Qell\).
  Let \((\Nbf_i)_i\) be a filtration of \(\Nbf\) by normal algebraic subgroups such that each \(\Nbf_i/\Nbf_{i+1}\) is a vector group centralized by \(\Nbf\).
  Let \(N\) be a compact open subgroup of \(\Nbf(\Qell)\) and \(N_i = N \cap \Nbf_i(\Qell)\).
  Then \((N, (N_i)_i)\) satisfies Assumption \ref{assu:N_alg}.
\end{lemm}
\begin{proof}
  Recall from \cite[\S IV.2.4]{DemazureGabriel} that unipotent algebraic groups over a field of characteristic zero are equivalent to nilpotent finite-dimensional Lie algebras via the exponential map and the Campbell-Hausdorff formula, which are polynomial in this case.
  Let \(\Bbf\) be the open affine subscheme of \(\Nbf_0^{r_0} \times_{\Qell} \dots \times_{\Qell} \Nbf_m^{r_m}\) parametrizing families \(\ul{g} = (g_{i,j})_{0 \leq i \leq m, 1 \leq j \leq r_i}\) such that for any \(0 \leq i \leq m\) the image of \((g_{i,j})_j\) in the vector group \(\Nbf_i/\Nbf_{i+1}\) is a basis.
  Then the morphism of schemes over \(\Qell\)
  \begin{align*}
    \Theta: \Bbf \times_{\Qell} \mathbb{A}_{\Qell}^{\sum_i r_i} & \longrightarrow \Bbf \times_{\Qell} \Nbf \\
    (\ul{g}, \ul{x}) & \longmapsto (\ul{g}, \psi_{\ul{g}}(\ul{x})),
  \end{align*}
  where \(\psi_{\ul{g}}\) is defined using the same formula \eqref{eq:def_psi_g}, is easily seen to be an isomorphism using the exponential map.
  Choose an arbitrary isomorphism \(\Nbf \simeq \mathbb{A}_{\Qell}^{\sum_i r_i}\) of schemes over \(\Qell\) (e.g.\ one given by a point in \(\Bbf(\Qell)\)), identifying \(\Theta\) with an automorphism of \(\mathbb{A}_{\Bbf}^{\sum_i r_i}\) over \(\Bbf\).
  Then the degrees of the specialization of \(\Theta\) and \(\Theta^{-1}\) at any point of \(\Bbf\) are uniformly bounded, implying the first bound in Assumption \ref{assu:N_alg}.

  The second bound is obtained using a similar argument for the composition of
  \begin{align*}
    \Bbf \times_{\Qell} \mathbb{A}_{\Qell}^{2\sum_i r_i} & \longrightarrow \Bbf \times_{\Qell} \Nbf \\
    (\ul{g}, \ul{y}, \ul{x}) & \longmapsto (\ul{g}, \psi_{\ul{g}}(\ul{y}) \psi_{\ul{g}}(\ul{x}))
  \end{align*}
  with \(\Theta^{-1}\).
\end{proof}

\begin{rema} \label{rem:assu_always_sat}
  One can show that any nilpotent, torsion-free and topologically finitely generated pro-\(\ell\) group \(N\) with a filtration \((N_i)_i\) as in Lemma \ref{lem:char_tors_free_nilp} \eqref{it:central_series_Zell_ri} satisfies Assumption \ref{assu:N_alg}.
  In fact any such pair \((N, (N_i)_i)\) is isomorphic to one obtained as in Lemma \ref{lem:alg_unip_satisfies_assu}.
  The proof of Mal'cev's theorem for discrete groups \cite[Theorems 17.2.5 and 17.3.2]{KargapolovMerzljakov_GTM} can certainly be adapted to this case.
  In the context of this paper it is natural to prove this result by induction using continuous cohomology groups in degree \(2\), see Corollary \ref{cor:all_N_sat_assu}.
\end{rema}

Under this assumption we may consider, for any \(\Ocal_E\)-module \(A\), the \(\Ocal_E\)-module \(\Pol(N^n, A)\) of polynomial functions \(N^n \to A\), defined as functions which are polynomial when composed with some \(\psi_{\ul{g}}\).
Because the polynomials \(\beta_i\) in \eqref{eq:change_basis_fil_N} are not affine in general it does not make sense to filter \(\Pol(N^n, A)\) by total degree.
To overcome this problem we introduce weights as follows.

\begin{defi} \label{def:weighted_nilp}
  A tuple \((N, (N_i)_i, \ul{w})\) will be called a filtered weighted nilpotent group if there exists \(d \in \Z\) such that
  \begin{itemize}
  \item \(N\) is a nilpotent, torsion-free and topologically generated pro-\(\ell\)-group,
  \item \((N_i)_{0 \leq i \leq m+1}\) is a filtration as in Lemma \ref{lem:char_tors_free_nilp} \eqref{it:central_series_Zell_ri},
  \item \((N, (N_i)_i)\) satisfies Assumption \ref{assu:N_alg} for the integer \(d\),
  \item \(\ul{w} = (w_A, w_m, \dots, w_0)\) are positive integers satisfying \(w_A \geq w_m\) and \(w_i \geq (d+1) w_{i-1}\) for all \(1 \leq i \leq m\).
  \end{itemize}
\end{defi}

Note that given \((N, (N_i)_i, d)\) satisfying Assumption \ref{assu:N_alg} there always exists \(\ul{w}\) satisfying the inequalities in Definition \ref{def:weighted_nilp}, for example we could take \(w_i = (d+1)^i\) and \(w_A = (d+1)^m\).

Let \((A, \Fil^\bullet A)\) be a filtered \(\Ocal_E\)-module.
Let \(\Fil^j_{\ul{w}} \Pol(N^n, A)\) be the sub-\(\Ocal_E\)-module of \(\Pol(N^n, A)\) generated by
\[ \left( \psi_{\ul{g}}(\ul{x}^{(1)}_m, \dots, \ul{x}^{(1)}_0), \dots, \psi_{\ul{g}}(\ul{x}^{(n)}_m, \dots, \ul{x}^{(n)}_0) \right) \mapsto P_m(\ul{x}^{(1)}_m, \dots, \ul{x}^{(n)}_m) \dots P_0(\ul{x}^{(1)}_0, \dots, \ul{x}^{(n)}_0) a \]
where \(P_i\) is a polynomial of total degree \(\leq d_i\) for some \(d_i \in \Z_{\geq 0}\), \(a \in \Fil^{d_A} A\) for some \(d_A \in \Z\), these upper bounds satisfying
\[ d_A w_A + d_m w_m + \dots + d_0 w_0 \leq j. \]
In other words \(x_{i,j}\) is given degree \(w_i\) and an element of \(\Fil^j A \smallsetminus \Fil^{j-1} A\) is given degree \(j w_A\).
In particular the ``error terms'' \(\beta_i\) in \eqref{eq:change_basis_fil_N} are now given smaller degree than the ``main terms'' \(\alpha_i\), and we conclude that the filtration \(\Fil^\bullet_{\ul{w}}\) on \(\Pol(N^n, A)\) does not depend on the choice of \(\ul{g}\).
This also shows that this filtration is preserved by any automorphism of the topological group \(N\) which preserves the filtration \((N_i)_i\).

A similar argument, now using the second part of Assumption \ref{assu:N_alg}, shows that if \((A, \Fil^\bullet A)\) is a filtered \(\Ocal_E[N]\)-module then the left action \eqref{eq:action_N_cochains} of \(N\) on \(\Pol(N^n, A)\) preserves \(\Fil^j_{\ul{w}} \Pol(N^n, A)\).
In this case we denote \(C^n_{\pol}(N, A) = \Pol(N^{n+1}, A)^N\) and endow it with the filtration \(\Fil^\bullet_{\ul{w}}\).

\begin{lemm} \label{lem:H_action_Fil_Pol}
  Let \((N, (N_i)_i, \ul{w})\) be a filtered weighted nilpotent group (Definition \ref{def:weighted_nilp}).
  Let \(H\) be a profinite topological group acting on \(N\) and preserving the filtration \((N_i)_i\).
  Let \((A, \Fil^\bullet A)\) be a filtered \(\Ocal_E[N \rtimes H]\)-module satisfying \(\Fil^j A = 0\) for some \(j \in \Z\).
  Let \(n \geq 0\).
  \begin{enumerate}
  \item Assume that the action of \(N\) on each \(\gr^j A\) is trivial.
    Then the action \eqref{eq:action_H_cochains} of \(H\) on \(\Pol(N^n, A)\) preserves the filtration \(\Fil^\bullet_{\ul{w}}\).
  \item Assume further that the action of \(H\) by conjugation on each \(N_i/N_{i+1}\) is trivial and that the action of \(N \rtimes H\) on each \(\gr^j A\) is trivial.
    Then \(H\) acts trivially on each \(\Fil^j_{\ul{w}} \Pol(N^n, A) / \Fil^{j-w_0}_{\ul{w}} \Pol(N^n, A)\).
  \end{enumerate}
\end{lemm}
\begin{proof}
  \begin{enumerate}
  \item For \(h \in H\) we know that \(\psi_{\ul{g}}^{-1}(h^{-1} \psi_{\ul{g}}(\ul{x}_m, \dots, \ul{x}_0) h)\) is given as in \eqref{eq:change_basis_fil_N} by a family \((\alpha_i)_i\) with \(\alpha_i \in \GL_{r_i}(\Zell)\) and a family of polynomials \((\beta_i)_i\) of degree at most \(d\) where \(d \in \Z_{\geq 0}\) is as in Definition \ref{def:weighted_nilp}.
    We have \(w_i \geq dw_j\) whenever \(i > j\) and the fact that \(H\) preserves each \(\Fil^j_{\ul{w}} \Pol(N^n, A)\) easily follows.
  \item Under the additional assumption we also have \(\alpha_i = 1\) for all \(i\) in \eqref{eq:change_basis_fil_N}.
    Using the inequalities \(w_i \geq dw_j + w_0\) whenever \(i>j\), as well as the inequality \(w_A \geq w_0\), we conclude that \(H\) acts trivially on each \(\Fil^j_{\ul{w}} \Pol(N^n, A) / \Fil^{j-w_0}_{\ul{w}} \Pol(N^n, A)\).
  \end{enumerate}
\end{proof}

We prove the higher-dimensional analogue of Lemma \ref{lem:Fil_Cpol_Fil_Pol_Zell}, but for a different purpose.

\begin{lemm} \label{lem:Fil_Cpol_Fil_Pol_N}
  Let \((N, (N_i)_i, \ul{w})\) be a filtered weighted nilpotent group (Definition \ref{def:weighted_nilp}).
  Let \((A, \Fil^\bullet A)\) be a filtered \(\Ocal_E[N]\)-module such that we have \(\Fil^j A = 0\) for some \(j \in \Z\) and each \(\Fil^j A\) is a finitely generated \(\Ocal_E\)-module with a continuous action of \(N\), and assume that the action of \(N\) on each graded pieces \(\gr^j A\) is trivial.
  Then for any \(d \in \Z\) and \(n \in \Z_{\geq 0}\) the map
  \begin{align*}
    \Fil^d_{\ul{w}} C^n_{\pol}(N, A) & \longrightarrow \Fil^d_{\ul{w}} \Pol(N^n, A) \\
    \varphi & \longmapsto ((x_1, \dots, x_n) \mapsto \varphi(1,x_1,\dots,x_n))
  \end{align*}
  is an isomorphism.
\end{lemm}
\begin{proof}
  We proceed by induction on \(\sum_i r_i \in \Z_{\geq 0}\), the initial case being trivial.
  So assume that we have \(\sum_i r_i > 0\).
  We may assume \(r_0>0\) and realize \(N\) as \(N' \rtimes \Zell\).
  By induction hypothesis we have isomorphisms
  \begin{align*}
    \Fil^d_{\ul{w}'} C^n_{\pol}(N', A) & \longrightarrow \Fil^d_{\ul{w}'} \Pol((N')^n, A) \\
    \varphi & \longmapsto ((x_1, \dots, x_n) \mapsto \varphi(1,x_1,\dots,x_n)).
  \end{align*}
  These isomorphisms are \(\Zell\)-equivariant (for the usual action given by formula \eqref{eq:action_H_cochains} on both sides), and we see them as an isomorphism of filtered \(\Ocal_E[\Zell]\)-modules
  \[ C^n_{\pol}(N', A) \longrightarrow \Pol((N')^n, A). \]
  By the second part of \ref{lem:H_action_Fil_Pol} we know that the action of \(\Zell\) on each \(\Fil^d_{\ul{w}'} / \Fil^{d-w_0}_{\ul{w}'}\) (on either side) is trivial.
  We also have identifications
  \begin{align*}
    \beta_n: \Pol(N^n, A) & \longrightarrow \Pol(\Zell^n, \Pol((N')^n, A)) \\
    \varphi & \longmapsto ((y_1, \dots, y_n) \mapsto ((x_1, \dots, x_n) \mapsto \varphi(x_1 y_1, \dots, x_n y_n)))
  \end{align*}
  which clearly identify the filtration \(\Fil^\bullet_{\ul{w}}\) on the left-hand side with \(\Fil^\bullet_{w_0}\) on the right-hand side (here \(\Pol((N')^n, A)\) is endowed with the filtration \(\Fil^\bullet_{\ul{w}'}\)).
  These maps are also \(N\)-equivariant (see Section \ref{sec:expl_HS}) and so we also have isomorphisms of filtered \(\Ocal_E\)-modules
  \begin{align*}
    \beta_{n+1}^N: C^n_{\pol}(N, A) & \longrightarrow C^n_{\pol}(\Zell, C^n_{\pol}(N', A)) \\
    \varphi & \longmapsto ((y_0, \dots, y_n) \mapsto ((x_0, \dots, x_n) \mapsto \varphi(x_0 y_0, \dots, x_n y_n))).
  \end{align*}
  We conclude the induction step using the commutative diagram of filtered \(\Ocal_E\)-modules
  \[
    \begin{tikzcd}
      C^n_{\pol}(N, A) \arrow[r, "{\beta_{n+1}^N}"] \arrow[dd, "{\alpha_N}"] & C^n_{\pol}(\Zell, C^n_{\pol}(N', A)) \arrow[d, "{C^n_{\pol}(\Zell, \alpha_{N'})}"] \\
      & \arrow[d, "{\alpha_{\Zell}}"] C^n_{\pol}(\Zell, \Pol((N')^n, A)) \\
      \Pol(N^n, A) \arrow[r, "{\beta_n}"] & \Pol(\Zell^n, \Pol((N')^n, A))
    \end{tikzcd}
  \]
  The maps \(\beta_n\) and \(\beta_{n+1}^N\) are isomorphisms of filtered \(\Ocal_E\)-modules, as are the maps \(C^n_{\pol}(\Zell, \alpha_{N'})\) (by induction hypothesis) and \(\alpha_{\Zell}\) (thanks to Lemma \ref{lem:Fil_Cpol_Fil_Pol_Zell}), so \(\alpha_N\) is also an isomorphism of filtered \(\Ocal_E\)-modules.
\end{proof}

\begin{theo} \label{thm:main_torsion}
  Let \((N, (N_i)_i, \ul{w})\) be a filtered weighted nilpotent group (Definition \ref{def:weighted_nilp}).
  Let \((A^\bullet, \Fil^\bullet A^\bullet)\) be a bounded below filtered complex of \(\Ocal_E[N]\)-modules satisfying the following conditions.
  \begin{itemize}
  \item There exists \(j \in \Z\) such that \(\Fil^j A^\bullet\) is the zero complex.
  \item For any \(j,n \in \Z\) the \(\Ocal_E\)-module \(\Fil^j A^n\) is finitely generated, the action of \(N\) on it is continuous, and the action of \(N\) on \(\gr^j A^n\) is trivial.
  \end{itemize}
  \begin{enumerate}
  \item Assume that each \(\Fil^j A^n\) is torsion, i.e.\ killed by some \(\mfrak_E^i\), and that there exists \(d_A \in \Z\) such that for any \(j \geq d_A\) the embedding \(\Fil^j A^\bullet \to \Fil^{j+1} A^\bullet\) is a quasi-isomorphism.
    Then for any integers \(d \geq d_A w_A + \sum_i r_i w_i\) and \(j \geq \lfloor d/w_A \rfloor\) the embedding of complexes
    \[ \Tot^\bullet(\Fil^d_{\ul{w}} C^\bullet_{\pol}(N, A^\bullet)) \to \Tot^\bullet(C^\bullet(N, \Fil^j A^\bullet)) \]
    is a quasi-isomorphism.
  \item Assume that there exists \(r_A \in \Z\) such that for any \(j \in \Z\) we have \(H^i(\Fil^j A^\bullet) = 0\) for \(i>r_A\).
    Then for any \(d \in \Z\) we have \(H^i(\Tot^\bullet(\Fil^d_{\ul{w}} C^\bullet_{\pol}(N, A^\bullet))) = 0\) for \(i>r_A+r_m+\dots+r_0\).
  \end{enumerate}
\end{theo}
\begin{proof}
  We will prove both points by induction on \(\sum_i r_i\).
  We start with the first point, in particular we assume that each \(\Fil^j A^n\) is torsion.
  We may reduce to the case where \(r_0>0\).
  (Note: this is the only place in the induction where we use \(w_A \geq w_m\) rather than the weaker inequality \(w_A \geq w_0\).)
  Let \(N'\) be the preimage of \(\Zell^{r_0-1} \times \{0\}\) in \(N\) and identify \(N\) with \(N' \rtimes \Zell\).
  By induction hypothesis for any \(d \geq d_A w_A + \sum_i r_i w_i - w_0\) and any \(j \geq \lfloor d/w_A \rfloor\) the embedding of complexes
  \[ \Tot^\bullet \left( \Fil^d_{\ul{w}'} C^\bullet_{\pol}(N', A^\bullet) \right) \rightarrow \Tot^\bullet \left( C^\bullet(N', \Fil^j A^\bullet) \right) \]
  is a quasi-isomorphism.
  Because for complexes of smooth \(\Ocal_E[\Zell]\)-modules the functor \(\Tot^\bullet(C^\bullet(\Zell, -))\) respect quasi-isomorphisms\footnote{This is a well-known and easier analogue of Lemma \ref{lem:derived_gp_coh_as_cocycles}.}, and by associativity of ``taking total complexes'', it follows that for such integers \(d\) and \(j\) the embedding
  \[ \Tot^\bullet \left( C^\bullet(\Zell, \Fil^d_{\ul{w}'} C^\bullet_{\pol}(N', A^\bullet)) \right) \rightarrow \Tot^\bullet \left( C^\bullet(\Zell, C^\bullet(N', \Fil^j A^\bullet)) \right) \]
  is also a quasi-isomorphism.
  Note that the induction hypothesis also implies that the embedding
  \[ \Tot^\bullet \left( \Fil^d_{\ul{w}'} C^\bullet_{\pol}(N', A^\bullet) \right) \longrightarrow \Tot^\bullet \left( \Fil^{d+1}_{\ul{w}'} C^\bullet_{\pol}(N', A^\bullet) \right) \]
  is a quasi-isomorphism for \(d \geq d_A w_A + \sum_i r_i w_i - w_0\).
  Thanks to Lemma \ref{lem:H_action_Fil_Pol} we know that the action \ref{eq:action_H_cochains} of \(\Zell\) on \(C^\bullet_{\pol}(N', A)\) preserves the filtration \(\Fil^j_{\ul{w}'}\) and that it acts trivially on each \(\Fil^j_{\ul{w}'} / \Fil^{j-w_0}_{\ul{w}'}\).
  By the first part of Lemma \ref{lem:Cpol_Zell_cplx} the embedding of complexes
  \[ \Tot^\bullet \left( \Fil^d_{w_0} C^\bullet_{\pol}(\Zell, C^\bullet_{\pol}(N', A^\bullet)) \right) \rightarrow \Tot^\bullet \left( C^\bullet(\Zell, \Fil^j_{\ul{w}'} C^\bullet_{\pol}(N', A^\bullet)) \right) \]
  is a quasi-isomorphism whenever \(d \geq d_A w_A + \sum_i r_i w_i\) and \(j \geq d\).

   By Theorem \ref{thm:Eilenberg_Zilber} the maps \(\AW\) and \(\EML\) induce quasi-isomorphisms
   \[ \Tot^\bullet( C^\bullet(\Zell, C^\bullet(N', A^n))) \leftrightarrow C^\bullet(N, A^n). \]
  We claim that for every \(n \in \Z\) the maps \(\AW\) and \(\EML\) also induce well defined maps
  \begin{equation} \label{eq:EZ_Fil}
    \Tot^\bullet( \Fil^d_{w_0} C^\bullet_{\pol}(\Zell, C^\bullet_{\pol}(N', A^n))) \leftrightarrow \Fil^d_{\ul{w}} C^\bullet_{\pol}(N, A^n)
  \end{equation}
  which are inverse of each other modulo homotopies, in particular they are quasi-isomorphisms.
   Here \(C^\bullet_{\pol}(N', A^n)\) is endowed with the filtration \(\Fil^\bullet_{\ul{w}'}\).
   To prove this claim we simply observe that it is clear on the defining formulas \eqref{eq:def_theta}, \eqref{eq:def_AW}, \eqref{eq:def_EML} and \eqref{eq:def_Psi} that the maps \(\theta^*\), \(\AW\), \(\EML\) and \(\Psi\) preserve the filtrations, so that the argument in Theorem \ref{thm:Eilenberg_Zilber} is still valid for these subcomplexes.

   We will conclude using the following commutative diagram where \(d\) is any integer and \(j \geq \lfloor d/w_A \rfloor\).
   \[
     \begin{tikzcd}
       \Tot^\bullet(\Fil^d_{w_0} C^\bullet_{\pol}(\Zell, C^\bullet_{\pol}(N', A^\bullet))) \arrow[r] \arrow[d, "{\AW}"] & \Tot^\bullet(C^\bullet(\Zell, \Fil^d_{\ul{w}'} C^\bullet_{\pol}(N', A^\bullet))) \arrow[d, "{\mathrm{IH}}"] \\
       \Tot^\bullet(\Fil^d_{\ul{w}} C^\bullet_{\pol}(N, A^\bullet)) \arrow[d, "{\iota}"] & \Tot^\bullet(C^\bullet(\Zell, \Tot^\bullet(C^\bullet(N', \Fil^j A^\bullet)))) \arrow[d, equal] \\
       \Tot^\bullet(C^\bullet(N, \Fil^j A^\bullet)) & \arrow[l, "{\AW}"] \Tot^\bullet(C^\bullet(\Zell, C^\bullet(N', A^\bullet)))
     \end{tikzcd}
   \]
   For \(d \geq d_A w_A + \sum_i r_i w_i\) we have seen above that the top horizontal map is a quasi-isomorphism thanks to Lemma \ref{lem:Cpol_Zell_cplx}, that the map labelled \(\mathrm{IH}\) is a quasi-isomorphism thanks to the induction hypothesis, and that both maps labelled \(\AW\) are quasi-isomorphisms.
   It follows that the embedding labelled \(\iota\) is also a quasi-isomorphism.

   The second point is proved by induction in a similar but simpler fashion using \eqref{eq:EZ_Fil} and the second part of Lemma \ref{lem:Cpol_Zell_cplx}.
\end{proof}

\begin{coro} \label{cor:derived_integral_main_result}
  Let \((N, (N_i)_i, \ul{w})\) be a filtered weighted nilpotent group (Definition \ref{def:weighted_nilp}).
  Let \(H\) be a profinite topological group acting on \(N\) and preserving \((N_i)_i\).
  Let \(A^\bullet\) be a bounded complex of objects of \(\Rep_{\fg,\cont}(N \rtimes H, \Ocal_E)\).
  Let \(\Fil^\bullet A^\bullet\) be a filtration of \(A^\bullet\) satisfying
  \begin{itemize}
  \item there exists \(j \in \Z\) for which we have \(\Fil^j A^\bullet = 0\),
  \item there exists \(d_A \in \Z\) for which \(\Fil^{d_A} A^\bullet = A^\bullet\),
  \item for any \(n \in \Z\) and \(j \in \Z\) the action of \(N\) on \(\gr^j A^n\) is trivial.
  \end{itemize}

  Then for any \(d \geq d_A w_A + \sum_i r_i w_i\) the obvious morphism in \(D^b(H, \Ocal_E)\)
  \[ F\left( \Tot^\bullet(\Fil^d_{\ul{w}} C^\bullet_{\pol}(N, A^\bullet)) \right) \to R\Gamma_{\cont}(N, F(A^\bullet)) \]
  is an isomorphism.
\end{coro}
\begin{proof}
    For \(n \in \Z\) and an integer \(i \geq 1\) we endow \(A^n / \mfrak_E^i A^n\) with the image of the given filtration of \(A^n\).
  By the first part of Theorem \ref{thm:main_torsion} we have for every \(i \geq 1\) a quasi-isomorphism
  \[ \Tot^\bullet \left( \Fil^d_{\ul{w}} C^\bullet_{\pol}(N, A^\bullet/\mfrak_E^i) \right) \to C^\bullet(N, A^\bullet/\mfrak_E^i) \]
  and as \(i\) varies these morphisms are clearly compatible.
  In particular we have an isomorphism in \(D^+(H, \Ocal_E)\)
  \[ \Tot^\bullet \left( \Fil^d_{\ul{w}} C^\bullet_{\pol}(N, F(A^\bullet)) \right) \to R\Gamma_{\cont}(N, F(A^\bullet)). \]
  Thanks to the second part of Theorem \ref{thm:main_torsion} we know that the complex \(\Tot^\bullet \left( \Fil^d_{\ul{w}} C^\bullet_{\pol}(N, A^\bullet) \right)\) has non-vanishing cohomology in only finitely many degrees, in particular we may apply to it the functor \(F\) of Proposition \ref{pro:F_to_Eke}.
  More precisely if we have \(A^j = 0\) for all \(j>r_A\) then we have
  \[ H^j \left( \Tot^\bullet \left( \Fil^d_{\ul{w}} C^\bullet_{\pol}(N, A^\bullet) \right) \right) =0 \]
  for all \(j > r_A + \sum_i r_i\), and letting \(r = 1 + r_A + \sum_i r_i\) we have a canonical isomorphism
  \[ F \left( \Tot^\bullet \left( \Fil^d_{\ul{w}} C^\bullet_{\pol}(N, A^\bullet) \right) \right) \simeq F \left( \tau_{\leq r} \Tot^\bullet \left( \Fil^d_{\ul{w}} C^\bullet_{\pol}(N, A^\bullet) \right) \right). \]
  (The integer \(r\) is one more than necessary here but this will be useful below.)
  We have a morphism in \(D^+(H, \Ocal_E)\)
  \begin{equation} \label{eq:F_FilCpol_to_FilCpol_F}
    F\left( \Tot^\bullet \left( \Fil^d_{\ul{w}} C^\bullet_{\pol}(N, A^\bullet) \right) \right) \to \Fil^d_{\ul{w}} C^\bullet_{\pol}(N, F(A^\bullet))
  \end{equation}
  obtained as the composition of the natural morphism of complexes
  \begin{equation} \label{eq:F_FilCpol_to_FilCpol_F_tau}
    F\left( \tau_{\leq r} \Tot^\bullet \left( \Fil^d_{\ul{w}} C^\bullet_{\pol}(N, A^\bullet) \right) \right) \to \tau_{\leq r} \Fil^d_{\ul{w}} C^\bullet_{\pol}(N, F(A^\bullet))
  \end{equation}
  with the quasi-isomorphism
  \[ \tau_{\leq r} \Fil^d_{\ul{w}} C^\bullet_{\pol}(N, F(A^\bullet)) \to \Fil^d_{\ul{w}} C^\bullet_{\pol}(N, F(A^\bullet)) \]
  (again using the second part of Theorem \ref{thm:main_torsion}, this time for \(A^\bullet \otimes_{\Ocal_E} \Ocal_E/\mfrak_E^i\)).
  Let us show that \eqref{eq:F_FilCpol_to_FilCpol_F} is an isomorphism in \(D^+(H, \Ocal_E)\).
  Thanks to Lemma \ref{lem:Fil_Cpol_Fil_Pol_N} we may identify the \(\Ocal_E\)-module \(\Fil^d_{\ul{w}} C^n_{\pol}(N,A^m)\) with a direct sum of finitely many \(\Fil^j A^m\).
  The natural maps of \((\Ocal_E)_\bullet\)-modules
  \[ \left( \Fil^j A^m \otimes_{\Ocal_E} \Ocal_E/\mfrak_E^i \right)_{i \geq 1} \longrightarrow \left( \left( \Fil^j A^m + \mfrak_E^i A^m \right) / \mfrak_E^i A^m \right)_{i \geq 1} \]
  are clearly surjective, and by the Artin-Rees lemma their kernels are essentially zero.
  This implies that the morphism of complexes \eqref{eq:F_FilCpol_to_FilCpol_F_tau} is surjective (in degree \(r\) this uses the fact that the right-hand side has vanishing cohomology, and is the reason for the increment in the definition of \(r\)) and that its kernel is essentially zero.
  In the derived category this kernel is isomorphic to the cone of \eqref{eq:F_FilCpol_to_FilCpol_F_tau} up to a shift, and so this cone is essentially zero and thus negligible by \cite[Proposition 3.4 (i)]{Ekedahl_adic}, i.e.\ \eqref{eq:F_FilCpol_to_FilCpol_F_tau} is an isomorphism in \(D^+(H, \Ocal_E)\) and thus so is \eqref{eq:F_FilCpol_to_FilCpol_F}.
\end{proof}

\begin{coro} \label{cor:cont_coh_unip_fg_OE_as_pol}
  Let \((N, (N_i)_i, \ul{w})\) and \((A^\bullet, \Fil^\bullet A^\bullet)\) be as in Theorem \ref{thm:main_torsion}.
  Assume that there exists \(d_A \in \Z\) for which we have \(\Fil^{d_A} A^\bullet = A^\bullet\).
  Then for any \(d \geq d_A w_A + \sum_i r_i w_i\) the embedding
  \[ \Tot^\bullet(\Fil^d_{\ul{w}} C^\bullet_{\pol}(N, A^\bullet)) \to \Tot^\bullet(C^\bullet_{\cont}(N, A^\bullet)) \]
  is a quasi-isomorphism, where on the right we consider continuous cochains for the \(\mfrak_E\)-adic topology on each \(A^n\).
\end{coro}
\begin{proof}
  By the usual double complexes arguments it is enough to consider the case of a single filtered \(\Ocal_E[N]\)-module \((A, \Fil^\bullet A)\).
  As in the proof of Theorem \ref{thm:main_torsion} we endow each \(A/\mfrak_E^iA\) with the filtration \(((\Fil^j A + \mfrak_E^i A)/\mfrak_E^i A)_{j \in \Z}\).
  From the first part of Theorem \ref{thm:main_torsion} we have quasi-isomorphisms
  \[ \Fil^d_{\ul{w}} C^\bullet_{\pol}(N, A/\mfrak_E^i A) \longrightarrow C^\bullet(N, A / \mfrak_E^i A) \]
  which are compatible as \(i \geq 1\) varies, and we consider them as an embedding of complexes of projective systems of abelian groups.
  All terms satisfy the Mittag-Leffler condition, in fact all transition morphisms are surjective (thanks to Lemma \ref{lem:Fil_Cpol_Fil_Pol_N} for the left-hand side), and so applying \(R \varprojlim_i\) we obtain a quasi-isomorphism
  \[ \varprojlim_i \Fil^d_{\ul{w}} C^\bullet_{\pol}(N, A/\mfrak_E^i A) \longrightarrow C^\bullet_{\cont}(N, A). \]
  Because each \(\Fil^j A\) is closed in \(A\) the natural map
  \[ \Fil^d_{\ul{w}} C^\bullet_{\pol}(N, A) \longrightarrow \varprojlim_i \Fil^d_{\ul{w}} C^\bullet_{\pol}(N, A/\mfrak_E^i A) \]
  is an isomorphism.
\end{proof}

\begin{coro} \label{cor:all_N_sat_assu}
  Let \(N\) be a pro-\(\ell\) group.
  Assume that \(N\) is topologically finitely generated, nilpotent and torsion-free.
  Let \((N_i)_{0 \leq i \leq m+1}\) be a filtration of \(N\) as in \ref{lem:char_tors_free_nilp} \eqref{it:central_series_Zell_ri}.
  Then there exists a unipotent algebraic group \(\Nbf\) over \(\Qell\) and a filtration \((\Nbf_i)_i\) of \(\Nbf\) by normal algebraic subgroups such that each \(\Nbf_i/\Nbf_{i+1}\) is a vector group centralized by \(\Nbf\), and an isomorphism between \(N\) and a compact open subgroup of \(\Nbf(\Qell)\) such that for any \(0 \leq i \leq m+1\) we have \(N_i = N \cap \Nbf_i(\Qell)\).
  In particular (Lemma \ref{lem:alg_unip_satisfies_assu}) Assumption \ref{assu:N_alg} is automatically satisfied.
\end{coro}
\begin{proof}
  We proceed by induction on \(m \geq 0\), the case \(m \leq 0\) being trivial.
  So assume \(m>0\) and consider \(N\) as a central extension of \(N' := N/N_m\) by \(N_m\).
  We have a continuous section \(s: N/N_m \to N\), for example the section given by Mal'cev coordinates \eqref{eq:def_psi_g}, and so our central extension corresponds to a class in \(H^2_{\cont}(N/N_m, N_m)\).
  Let \(N'_i = N_i/N_m\) for \(0 \leq i \leq m\), applying the induction hypothesis to \((N', (N'_i)_i)\) we obtain a filtered unipotent group \((\Nbf', (\Nbf'_i)_{0 \leq i \leq m})\).
  By Lemma \ref{lem:alg_unip_satisfies_assu} there exists \(\ul{w}' = (w_A, w_{m-1}, \dots, w_0)\) such that \((N', (N'_i)_i, \ul{w}')\) is a weighted filtered nilpotent group (Definition \ref{def:weighted_nilp}).
  By Corollary \ref{cor:cont_coh_unip_fg_OE_as_pol} any class in \(H^2_{\cont}(N/N_m, N_m)\) is represented by a polynomial \(2\)-cocycle, i.e.\ up to choosing another continuous section \(s\) we may assume that the map
  \begin{align*}
    (N/N_m)^2 & \longmapsto N_m \\
    (x,y) & \longmapsto s(x) s(y) s(xy){-1}
  \end{align*}
  is given by a polynomial, and choosing an isomorphism \(N_m \simeq \Zell^{r_m}\) we can see this polynomial as a morphism \(c: (\Nbf')^2 \longrightarrow \mathbb{A}_{\Qell}^{r_m}\) of schemes over \(\Qell\).
  Using \(c\) we may define a structure of linear algebraic group on the scheme \(\Nbf = \mathbb{A}_{\Qell}^{r_m} \times_{\Qell} \Nbf'\), define \(\Nbf_m = \mathbb{A}_{\Qell}^{r_m} \subset \Nbf\) and \(\Nbf_i = \mathbb{A}_{\Qell}^{r_m} \times_{\Qell} \Nbf'_i\) for \(0 \leq i < m\) and identify \(N\) with a compact open subgroup of \(\Nbf(\Qell)\).
  Details are left to the reader.
\end{proof}

\begin{rema}
  In the same vein using cohomology in degree \(1\) one can show that any continuous unipotent finite-dimensional representation over \(E\) of a nilpotent, torsion-free and topologically finitely generated pro-\(\ell\) group \(N\) (i.e.\ any compact open subgroup of \(\Nbf(\Qell)\) for a unipotent linear algebraic group \(\Nbf\) over \(\Qell\)) is automatically algebraic, by induction on the dimension.
\end{rema}

\section{Characteristic zero coefficients}
\label{sec:char_0_coeff}

Let \(\Nbf\) be a unipotent linear algebraic group over \(\Qell\).
Let \((\Nbf_i)_{0 \leq i \leq m+1}\) be a filtration by normal algebraic subgroups such that we have \(\Nbf_0 = \Nbf\) and \(\Nbf_{m+1} = 1\) and such that each \(\Nbf_i/\Nbf_{i+1}\) is a vector group centralized by \(\Nbf\).
Recall from the proof of Lemma \ref{lem:alg_unip_satisfies_assu} that we have analogues of Mal'cev coordinates \ref{eq:def_psi_g} for \((\Nbf(\Qell), (\Nbf_i(\Qell))_i)\) with \(\Qell\) replacing \(\Zell\), and that they satisfy the analogue of Assumption \ref{assu:N_alg}.
In particular there exists \(\ul{w} = (w_A, w_m, \dots, w_0)\) as in Definition \ref{def:weighted_nilp} but now for \((\Nbf, (\Nbf_i)_i)\).
Similarly to Definition \ref{def:weighted_nilp} we will call such a triple \((\Nbf, (\Nbf_i)_{0 \leq i \leq m+1}, \ul{w})\) a weighted filtered unipotent linear algebraic group.
Let \(V\) be a (finite-dimensional) algebraic representation of \(\Nbf_E\), which we consider as a continuous representation of \(\Nbf(\Qell)\).
Then \(V\) admits a filtration \(\Fil^\bullet V\) by subrepresentations of \(\Nbf_E\) such that there exists \(j \in \Z_{>0}\) satisfying \(\Fil^{-j} V = 0\) and \(\Fil^j V = V\) and \(\Nbf_E\) acts trivially on each \(\gr^j V\).
We may define \(\Pol(\Nbf(\Qell)^n, V)\) similarly to the integral case after Definition \ref{def:weighted_nilp}, endow it with a filtration \(\Fil^\bullet_{\ul{w}}\).
Similarly we can define the \(E\)-vector space \(C^n_{\pol}(\Nbf(\Qell), V)\) with its filtration \(\Fil^\bullet_{\ul{w}}\).
For any compact open subgroup \(N\) of \(\Nbf(\Qell)\) there exists an \(\Ocal_E\)-lattice \(A\) in \(V\) which is stable under \(N\), and we endow it with the induced filtration.
Then each \(\Fil^j_{\ul{w}} \Pol(N^n, A)\) (resp.\ \(\Fil^j_{\ul{w}} C^n_{\pol}(N, A)\)) is an \(\Ocal_E\)-lattice in \(\Fil^j_{\ul{w}} \Pol(\Nbf(\Qell)^n, V)\) (resp.\ \(\Fil^j_{\ul{w}} C^n_{\pol}(\Nbf(\Qell), V)\)), essentially because \(N\) is Zariski-dense in \(\Nbf\).
We can use this simple fact, combined with the observation that compact open subgroups of \(\Nbf(\Qell)\) exist, to deduce analogues of some of the results of the previous section with \(\Nbf(\Qell)\) replacing \(N\) and algebraic representations of \(\Nbf_E\) replacing finitely generated \(\Ocal_E\)-modules with a continuous action of \(N\).
In particular Lemmas \ref{lem:H_action_Fil_Pol} and \ref{lem:Fil_Cpol_Fil_Pol_N} admit obvious analogues, as does the second part of Theorem \ref{thm:main_torsion}.
We record a slightly indirect consequence of our integral results in the following lemma.

\begin{lemm} \label{lem:alg_N_Fil_Cpol_stabilizes}
  Let \((\Nbf, (\Nbf_i)_{0 \leq i \leq m+1})\) be as above.
  Let \((V^\bullet, \Fil^\bullet V^\bullet)\) be a bounded below filtered complex of finite-dimensional algebraic representations of \(\Nbf_E\) such that there exists \(j \in \Z\) satisfying \(\Fil^j V^\bullet = 0\) and such that \(\Nbf_E\) acts trivially on each graded piece \(\gr^j V^n\).
  Assume that there exists \(d_V \in \Z\) satisfying \(\Fil^{d_V} V^\bullet = V^\bullet\).
  Then for any \(j \geq d_V w_A + \sum_i r_i w_i\) the embedding
  \[ \Tot^\bullet(\Fil^j_{\ul{w}} C^\bullet_{\pol}(\Nbf(\Qell), V^\bullet)) \longrightarrow \Tot^\bullet(\Fil^{j+1}_{\ul{w}} C^\bullet_{\pol}(\Nbf(\Qell), V^\bullet)) \]
  is a quasi-isomorphism.
\end{lemm}
\begin{proof}
  Choose a compact open subgroup \(N\) of \(\Nbf(\Qell)\) and let \(N_i = N \cap \Nbf_i(\Qell)\).
  Choose \(\Ocal_E\)-lattices \(A^n\) in \(V^n\) which are stable under \(N\) and stable under the differentials of \(V^\bullet\), and endow \(A^\bullet\) with the induced filtration.
  The lemma follows from Corollary \ref{cor:cont_coh_unip_fg_OE_as_pol} and extending scalars from \(\Ocal_E\) to \(E\).
\end{proof}

The following lemma will be useful to apply Corollary \ref{cor:derived_integral_main_result} in cases where instead of a semi-direct product \(N \rtimes H\) we have an open subgroup of such a product, provided we invert \(\ell\) in the coefficients.

\begin{lemm} \label{lem:res_iso_after_inv_ell}
  Let \(K'\) be an open subgroup of a profinite topological group \(K\).
  Let \(N\) be a closed normal subgroup of \(K\).
  Assume that \(N\) is pro-\(\ell\), nilpotent, torsion-free and topologically finitely generated (see Lemma \ref{lem:char_tors_free_nilp}).
  Denote \(N' = K' \cap N\), a closed normal subgroup of \(K'\) and an open subgroup of \(N\).
  Let \(A^\bullet\) be a bounded complex of objects of \(\Rep_{\fg,\cont}(K, \Ocal_E)\).
  Assume the existence of a filtration \((\Fil^j A^\bullet)_{j \in \Z}\) of \(A^\bullet\) satisfying
  \begin{itemize}
  \item for any \(n \in \Z\) and any \(j \in \Z\) the action of \(N\) on \(\gr^j A^n\) is trivial,
  \item there exists \(j \in \Z_{\geq 0}\) satisfying \(\Fil^{-j} A^\bullet = 0\) and \(\Fil^j A^\bullet = A^\bullet\).
  \end{itemize}
  Then the morphism (see \eqref{eq:res_RGamma})
  \[ r_{K,N,K'}(F(A^\bullet)): \res_{K'/N'} R\Gamma(N, F(A^\bullet)) \rightarrow R\Gamma(N', \res_{K'} F(A^\bullet)) \]
  becomes an isomorphism in \(D^+(S_{K'/N'}^{\N}, (\Ocal_E)_\bullet)[\ell^{-1}]\), in particular it becomes an isomorphism in \(D^+(K'/N', E)\).
\end{lemm}
\begin{proof}
  By the same argument as in the proof of Lemma \ref{lem:triang_cat_inv_ell} it is enough to show that the cone of \(r_{K,N,K'}(F(A^\bullet))\) has non-zero cohomology in only finitely many degrees and that all of its cohomology group are killed by some power of \(\ell\).
  For this we can forget the action of \(K'/N'\), and using Lemma \ref{lem:cmp_res_inj_cocyc} we are left to show this property for the cone of the restriction map on cocycles
  \[ R\Gamma_{\cont}(N, \res_N F(A^\bullet)) \longrightarrow R\Gamma_{\cont}(N', \res_{N'} F(A^\bullet)). \]
  Let \((N_i)_i\) be the upper central series of \(N\), and choose \(\ul{w}\) as in Definition \ref{def:weighted_nilp}.
  As in the proof of Corollary \ref{cor:derived_integral_main_result} we endow each \(A^n/\mfrak_E^i A^n\) with the filtration induced by that on \(A^n\), and for any large enough integer \(d\) we have thanks to the first part of Theorem \ref{thm:main_torsion} a quasi-isomorphism
  \[ \Tot^\bullet(\Fil^d_{\ul{w}} C^\bullet_{\pol}(N, F(A^\bullet))) \simeq R\Gamma_{\cont}(N, \res_N F(A^\bullet)), \]
  and similarly with \(N\) replaced by \(N'\) (replacing \(N_i\) by \(N_i' := N_i \cap N'\)).
  We are left to consider the cone of the natural map of complexes
  \begin{equation} \label{eq:res_cocyc_cone}
    \Tot^\bullet(\Fil^d_{\ul{w}} C^\bullet_{\pol}(N, F(A^\bullet))) \longrightarrow \Tot^\bullet(\Fil^d_{\ul{w}} C^\bullet_{\pol}(N', F(A^\bullet))).
  \end{equation}
  By the second part of Theorem \ref{thm:main_torsion} we also know that both sides have non-vanishing cohomology in only finitely many degrees.

  We claim that for any \(m \geq 0\) and \(n \in \Z\) there exists an integer \(k \geq 0\), which also depends on \(N\), \(N'\), \(\ul{w}\) and on an integer \(j_0>0\) satisfying \(\Fil^{-j_0} A^n = 0\) and \(\Fil^{j_0} A^n = A^n\), such that for any \(i \geq 1\) the kernel and cokernel of the morphism
  \[ \Fil^d_{\ul{w}} C^m_{\pol}(N, A^n/\mfrak_E^i A^n) \xrightarrow{\res} \Fil^d_{\ul{w}} C^m_{\pol}(N', A^n/\mfrak_E^i A^n) \]
  are killed by \(\ell^k\).
  To prove the claim we use Lemma \ref{lem:Fil_Cpol_Fil_Pol_N} to consider instead the kernel and cokernel of
  \begin{equation} \label{eq:res_Fil_Pol_N_Nprime}
    \Fil^d_{\ul{w}} \Pol(N^m, A^n/\mfrak_E^i A^n) \xrightarrow{\res} \Fil^d_{\ul{w}} \Pol((N')^m, A^n/\mfrak_E^i A^n).
  \end{equation}
  We also observe that we have an isomorphism of filtered \(\Ocal_E\)-modules
  \[ \Pol(N^m, \Ocal_E) \otimes_{\Ocal_E} A^n/\mfrak_E^i A^n \simeq \Pol(N^m, A^n/\mfrak_E^i A^n), \]
  where \(\Ocal_E\) is endowed with the trivial filtration (\(\Fil^{-1} \Ocal_E = 0\) and \(\Fil^0 \Ocal_E = \Ocal_E\)) and the filtration on \(A^n/\mfrak_E^i A^n\) is divided by \(w_A\), and similarly for \(N'\).
  Each restriction map \(\Fil^j_{\ul{w}} \Pol(N^m, \Ocal_E) \to \Fil^j_{\ul{w}} \Pol((N')^m, \Ocal_E)\) is an embedding with finite cokernel because both sides are lattices inside \(\Fil^j_{\ul{w}} \Pol(\Nbf(\Qell)^m, E)\).
  It follows that the cokernel of \eqref{eq:res_Fil_Pol_N_Nprime} is killed by a power of \(\ell\) which does not depend on \(i\) (note that only the indices \(j\) satisfying \(0 \leq j \leq d + (j_0-1) w_A\) are relevant).
  The fact that the kernel of \eqref{eq:res_Fil_Pol_N_Nprime} is killed by some power of \(\ell\) is slightly more subtle, it follows by induction from the same property for the maps
  \begin{equation} \label{eq:res_gr_Pol_N_Nprime}
    \gr^j_{\ul{w}} \Pol(N^m, A^n/\mfrak_E^i A^n) \xrightarrow{\res} \gr^j_{\ul{w}} \Pol((N')^m, A^n/\mfrak_E^i A^n)
  \end{equation}
  for \(w_A(1-j_0) \leq j \leq d\).
  We have
  \[ \gr^j_{\ul{w}} \Pol(N^m, A^n/\mfrak_E^i A^n) \simeq \bigoplus_{\substack{j_P,j_A \in \Z \\ j_P + w_A j_A = j}} \gr^{j_P}_{\ul{w}} \Pol(N^m, \Ocal_E) \otimes_{\Ocal_E} \gr^{j_A} A^n/\mfrak_E^i A^n \]
  and similarly for \(N'\).
  Using the beginning of the long exact sequence for \(\mathrm{Tor}\) groups we deduce that if a power of \(\ell\) kills
  \[ \gr^{j_P}_{\ul{w}} \Pol((N')^m, \Ocal_E) / \gr^{j_P}_{\ul{w}} \Pol(N^m, \Ocal_E) \]
  for all \(0 \leq j_P \leq j+(j_0-1)w_A\) then it also kills the kernel of \eqref{eq:res_gr_Pol_N_Nprime}.
  This concludes the proof of the claim.
  It is easy to deduce that each cohomology group of the cone of \eqref{eq:res_cocyc_cone} is killed by some power of \(\ell\).
\end{proof}

We can now state and prove a slight generalization of Corollary \ref{cor:derived_integral_main_result} when \(\ell\) is inverted in the coefficients.

\begin{prop} \label{pro:cont_coh_as_pol_ell_inv}
  Let \((\Nbf, (\Nbf_i)_{0 \leq i \leq m+1}, \ul{w})\) be a filtered weighted unipotent linear algebraic group over \(\Qell\).
  Let \(\Hbf\) be a linear algebraic group over \(\Qell\) acting on \(\Nbf\).
  Let \(V^\bullet\) be a bounded complex of finite-dimensional algebraic representations of \((\Nbf \rtimes \Hbf)_E\), considered as continuous \(E\)-representations of \((\Nbf \rtimes \Hbf)(\Qell)\).
  Let \(\Fil^\bullet V^\bullet\) be a filtration on \(V^\bullet\) satisfying
  \begin{itemize}
  \item there exists \(j \in \Z\) for which we have \(\Fil^j V^\bullet = 0\),
  \item there exists \(d_V \in \Z\) for which we have \(\Fil^{d_V} V^\bullet = V^\bullet\),
  \item for any \(n \in \Z\) and \(j \in \Z\) the action of \(\Nbf_E\) on \(\gr^j V^n\) is trivial.
  \end{itemize}
  
  Let \(K\) be a compact open subgroup of \((\Nbf \rtimes \Hbf)(\Qell)\).
  Denote \(N_K = \Nbf(\Qell) \cap K\) and \(H_K = K/N_K\), considered as a compact open subgroup of \(\Hbf(\Qell)\).

  Then for any \(d \geq d_V w_A + \sum_i r_i w_i\) we have an isomorphism in \(D^+(H_K, E)\)
  \[ R\Gamma(N_K, F(V^\bullet)) \simeq F(\Tot^\bullet(\Fil^d_{\ul{w}} C^\bullet_{\pol}(\Nbf(\Qell), V^\bullet))), \]
  where the right-hand side is well-defined because each \(\Fil^d_{\ul{w}} C^\bullet_{\pol}(\Nbf(\Qell), V^n)\) is exact in degree greater than \(\sum_i r_i\).
\end{prop}
Note that for any complex \(V^\bullet\) as in the proposition there always exists a filtration satisfying the assumptions in the proposition because every algebraic representation of \(\Nbf\) is unipotent.
\begin{proof}
  Let \(s: \Hbf \to \Nbf \rtimes \Hbf\) be the obvious section.
  We introduce the compact open subgroup \(K' = K s(H_K)\) of \((\Nbf \rtimes \Hbf)(\Qell)\).
  It is open because it contains \(K\), and compact because it is generated by its open subgroup \(K\) and the subset \(X := \{k s(k)^{-1} | k \in K \}\) of \(\Nbf(\Qell)\), and \(X\) is compact and so there exists a compact subgroup of \(\Nbf(\Qell)\) which contains \(X \cup N_K\) and is stable under the conjugation action of \(K\).
  Denoting \(N_{K'} = \Nbf(\Qell) \cap K'\), we have \(K' = N_{K'} \rtimes s(H_K)\), in particular \(H_{K'} := K'/N_{K'}\) is equal to \(H_K\).

  Let \(A^\bullet\) be a model of \(V^\bullet\), i.e.\ an object of \(D^b(\Rep_{\fg,\cont}(K', \Ocal_E))\) with an isomorphism \(E \otimes_{\Ocal_E} A^\bullet \simeq V^\bullet\) in \(D^b(\Rep_{\fg,\cont}(K', E))\).
  We may and do assume that \(A^\bullet\) is simply a sub-complex of \(V^\bullet\) considered as a complex of \(\Ocal_E[K']\)-modules.
  Define \(\Fil^j A^n = A^n \cap \Fil^j V^n\) for all \(n,j \in \Z\).
  By Lemma \ref{lem:derived_gp_coh_as_cocycles} and Corollary \ref{cor:derived_integral_main_result} we have an isomorphism in \(D^+(H_K, \Ocal_E)\)
  \[ R\Gamma(N_{K'}, F(A^\bullet)) \simeq F \left( \Tot^\bullet(\Fil^d_{\ul{w}} C^\bullet_{\pol}(N_{K'}, A^\bullet)) \right). \]
  In particular this becomes an isomorphism in \(D^+(H_K, E)\).
  It is clear that the image of
  \[ \Tot^\bullet(\Fil^d_{\ul{w}} C^\bullet_{\pol}(N_{K'}, A^\bullet)) \]
  in \(D^b(\Rep_{\fg,\cont}(H_K, E))\) is equal to
  \[ \Tot^\bullet(\Fil^d_{\ul{w}} C^\bullet_{\pol}(\Nbf(\Qell), V^\bullet)). \]
  Finally by Lemma \ref{lem:res_iso_after_inv_ell} we have an isomorphism in \(D^+(H_K, E)\)
  \[ R\Gamma(N_{K'}, F(A^\bullet)) \simeq R\Gamma(N_K, \res_K F(A^\bullet)). \]
\end{proof}

\section{Lie algebra cohomology}
\label{sec:Lie_alg_coh}

In order to work in derived categories of finite-dimensional vector spaces over \(E\) with additional structure, we define Lie algebra cohomology in an ad hoc manner using the Chevalley-Eilenberg complex as follows.
Let \(\Nbf\) be a linear algebraic group over \(E\), and let \(\nfrak\) be its Lie algebra.
Let \(\Hbf\) be another linear algebraic group over \(E\).
Let \(\Rep_{\fg, \alg}(\Nbf \rtimes \Hbf)\) be the abelian category of finite-dimensional vector spaces over \(E\) endowed with an algebraic linear action of \(\Nbf \rtimes \Hbf\).
For an object \(V\) in this category, recall the Chevalley-Eilenberg complex of \(E\)-vector spaces \(C^\bullet(\nfrak, V) = \Hom_E(\bigwedge^\bullet \nfrak, V)\) with differential
\begin{multline*}
  df(X_1, \dots, X_{n+1}) = \sum_{i=1}^{n+1} (-1)^{i+1} X_i \cdot f(X_1, \dots, \widehat{X_i}, \dots, X_{n+1}) \\
  + \sum_{1 \leq i < j \leq n+1} (-1)^{i+j} f([X_i, X_j], X_1, \dots, \widehat{X_i}, \dots, \widehat{X_j}, \dots, X_{n+1}).
\end{multline*}
The action of \(\Hbf\) on \(C^\bullet(\nfrak, V)\) is defined by the formula
\begin{equation} \label{eq:action_on_CE}
 (y \cdot f)(X_1, \dots, X_n) = y \cdot f(y^{-1} X_1 y, \dots, y^{-1} X_n y)
\end{equation}
where \(y \in \Hbf(R)\) for some commutative \(E\)-algebra \(R\), \(X_1, \dots, X_n \in \nfrak \otimes_E R\) and \(f \in R \otimes_E C^n(\nfrak, V)\).
Define
\begin{align*}
  R\Gamma_{\Lie}(\nfrak, -): D^b(\Rep_{\fg, \alg}(\Nbf \rtimes \Hbf)) & \longrightarrow D^b(\Rep_{\fg, \alg}(\Hbf)) \\
  V^\bullet & \longmapsto \Tot^\bullet \left( \Hom \left( \bigwedge^\bullet \nfrak, V^\bullet \right) \right).
\end{align*}

The fact that this is well-defined, i.e.\ that this functor preserves quasi-isomorphisms can be checked using the commutative diagram
\[
  \begin{tikzcd} [column sep=6em]
    K^b(\Rep_{\fg,\alg}(\Nbf \rtimes \Hbf)) \arrow[r, "{R\Gamma_{\Lie}(\nfrak, -)}"] \arrow[d] & K^b(\Rep_{\fg,\alg}(\Hbf)) \arrow[d] \\
    K^+(U(\nfrak)\mathrm{-Mod}) \arrow[r, "{R\Gamma_{\Lie}(\nfrak, -)}"] & K^+(E\mathrm{-Vect})
  \end{tikzcd}
\]
where the two vertical (forgetful) functors become conservative after passing to derived categories and the bottom horizontal functor is known to compute Lie algebra cohomology as a derived functor.

The analogue of Theorem \ref{thm:Eilenberg_Zilber} in the setting of Lie algebras is much simpler.

\begin{lemm} \label{lem:HS_Lie_adhoc}
  Consider a compound semi-direct product \(N_1 \rtimes (N_2 \rtimes H)\) of linear algebraic groups over \(E\).
  Let \(N = N_1 \rtimes N_2\).
  For any object \(V\) of \(\Rep_{\fg, \alg}(N_1 \rtimes (N_2 \rtimes H))\) we have an isomorphism of complexes
  \begin{align*}
    C^\bullet(\nfrak, V) & \longrightarrow \Tot^\bullet(C^\bullet(\nfrak_2, C^\bullet(\nfrak_1, V))) \\
    f & \longmapsto \left( (Y_1, \dots, Y_i) \mapsto ((X_1, \dots, X_j) \mapsto f(X_1, \dots, X_j, Y_1, \dots, Y_i)) \right).
  \end{align*}

  In particular this gives an isomorphism between the two functors \(R\Gamma_{\Lie}(\nfrak, -)\) and \(R\Gamma_{\Lie}(\nfrak_2, R\Gamma_{\Lie}(\nfrak_1, -))\)
  \[ D^b(\Rep_{\fg, \alg}(\Nbf_1 \rtimes (\Nbf_2 \rtimes \Hbf))) \longrightarrow D^b(\Rep_{\fg, \alg}(\Hbf)). \]
\end{lemm}
\begin{proof}
  We have \(\nfrak = \nfrak_1 \oplus \nfrak_2\) with \([\nfrak_k, \nfrak_k] \subset \nfrak_k\) and \([\nfrak_1, \nfrak_2] \subset \nfrak_1\).
  The action of \(\nfrak_2\) on \(C^\bullet(\nfrak_1, V)\) is easily computed from
  \eqref{eq:action_on_CE}:
  \[ (Y \cdot f)(X_1, \dots, X_n) = Y \cdot f(X_1, \dots, X_n) + \sum_{i=1}^n f(X_1, \dots, [X_i,Y], \dots, X_n). \]
  We have isomorphisms
  \begin{align*}
    \bigoplus_{i+j=n} \bigwedge^j \nfrak_1 \otimes \bigwedge^i \nfrak_2 & \longrightarrow \bigwedge^n \nfrak \\
    (X_1 \wedge \dots \wedge X_j) \otimes (Y_1 \wedge \dots \wedge Y_i) & \longmapsto X_1 \wedge \dots \wedge X_j \wedge Y_1 \wedge \dots \wedge Y_i
  \end{align*}
  and so the map in the lemma is an isomorphism of rational representations of \(\Hbf\).
  To check that it is a map of complexes, if \(f \in C^{n-1}(\nfrak, V)\) corresponds to \((f_{i,j})_{i+j=n-1}\) then for all \(X_1, \dots, X_j \in \nfrak_1\) and \(Y_1, \dots, Y_i \in \nfrak_2\) we have
  \[ f(X_1, \dots, X_j, Y_1, \dots, Y_i) = f_{i,j}(Y_1, \dots, Y_i)(X_1, \dots, X_j) \]
  and we compute
  \begin{align*}
    & (d_{\nfrak} f)(X_1, \dots, X_j, Y_1, \dots, Y_i) \\
    =& \sum_{a=1}^j (-1)^{a+1} X_a \cdot f_{i,j-1}(Y_1, \dots, Y_i)(X_1, \dots, \widehat{X_a}, \dots, X_j) \\
    &+ \sum_{a=1}^i (-1)^{j+a+1} Y_a \cdot f_{i-1,j}(Y_1, \dots, \widehat{Y_a}, \dots, Y_i)(X_1, \dots, X_j) \\
    &+ \sum_{1 \leq a < b \leq j} (-1)^{a+b} f_{i,j-1}(Y_1, \dots, Y_i)([X_a,X_b], X_1, \dots, \widehat{X_a}, \dots, \widehat{X_b}, \dots, X_j) \\
    &+ \sum_{1 \leq a < b \leq i} (-1)^{a+b+j} f_{i-1,j}([Y_a,Y_b], Y_1, \dots, \widehat{Y_a}, \dots, \widehat{Y_b}, \dots, Y_i)(X_1, \dots, X_j) \\
    &+ \sum_{1 \leq a \leq j} \sum_{1 \leq b \leq i} (-1)^{a+b+j} f_{i-1,j}(Y_1, \dots, \widehat{Y_b}, \dots, Y_i)([X_a,Y_b], X_1, \dots, \widehat{X_a}, \dots, X_j) \\
    =& d_{\nfrak_1} (f_{i,j-1}(Y_1, \dots, Y_i))(X_1, \dots, X_j) + (-1)^j (d_{\nfrak_2} f_{i-1,j})(Y_1, \dots, Y_i)(X_1, \dots, X_j).
  \end{align*}
  The second part of the lemma follows using standard arguments with total complexes.
\end{proof}

\section{Algebraic group cohomology and Lie algebra cohomology}
\label{sec:alg_gp_coh_and_Lie}

Let \(\Nbf\) be a connected linear algebraic group over \(E\).
Denote \(\Nbf = \Spec R(\Nbf)\).
Following Hochschild \cite{Hochschild_cohalglingps} we consider the abelian category \(\Rep_{\alg}(\Nbf)\) of representations of \(N\) which are sums (i.e.\ internal direct limits) of finite-dimensional algebraic representations.
This category has enough injectives: for any object \(V\) we have an embedding \(V \to R(\Nbf) \otimes_E V\), and \(R(\Nbf) \otimes_E V\) is injective in this category because for any object \(W\) we have an identification
\begin{align*}
  \Hom_{\Nbf}(W, R(\Nbf) \otimes_E V) & \longrightarrow \Hom_E(W, V) \\
  f & \longmapsto (w \mapsto f(w)(1))
\end{align*}
essentially because the equivariance condition reads \(f(w)(g) = g \cdot f(g^{-1} \cdot w)(1)\).
Just as in the case of group cohomology we have a resolution of an arbitrary object \(V\) of \(\Rep_{\alg}(\Nbf)\):
\[ V \to R(\Nbf) \otimes_E V \to R(\Nbf^2) \otimes_E V \to \dots \]
with differentials as in \eqref{eq:diff_hom_cochains}, that is
\[ d(f_0 \otimes \dots \otimes f_n \otimes a) = \sum_{j=0}^{n+1} (-1)^j f_0 \otimes \dots \otimes f_{j-1} \otimes 1 \otimes f_j \otimes \dots \otimes f_n \otimes a. \]
Noticing that each term of this resolution is an injective object because of the identification \(R(\Nbf^{i+1}) \simeq R(\Nbf) \otimes_E R(\Nbf^i)\), we conclude that the right derived functor of the functor of invariants applied to the complex \(V\) concentrated in degree \(0\) is (up to canonical quasi-isomorphism) the complex
\[ C^\bullet_{\alg}(\Nbf, V) := (R(\Nbf^{\bullet+1}) \otimes_E V)^{\Nbf} \]
supported in non-negative degrees.
Note that if \(\Nbf = \Ubf_E\) for some connected linear algebraic group \(\Ubf\) over \(\Qell\) and if \(V\) is finite-dimensional then we have a canonical isomorphism of complexes
\[ C^\bullet_{\alg}(\Nbf, V) \simeq C^\bullet_{\pol}(\Ubf(\Qell), V) \]
where the right-hand side is the complex introduced in Section \ref{sec:char_0_coeff}.
This follows from the well-known fact that \(\Ubf(\Qell)\) is Zariski-dense in \(\Ubf\).

There is another natural complex to which an arbitrary object \(V\) of \(\Rep_{\alg}(\Nbf)\) maps.
Denote by \(\Omega^\bullet_{\Nbf/E}\) the algebraic de Rham complex on \(\Nbf\).
There is a natural action of \(\Nbf\) on each term by left translation, commuting with the differentials, and it is easy to check that each term is an object of \(\Rep_{\alg}(\Nbf)\).
Now for any object \(V\) of \(\Rep_{\alg}(\Nbf)\) we obtain a complex \(\Omega^\bullet_{\Nbf/E} \otimes_E V\) and the kernel of \(\Omega^0_{\Nbf/E} \otimes_E V \to \Omega^1_{\Nbf/E} \otimes_E V\) is identified with \(V\).
It is well-known that the multiplication map
\[ R(\Nbf) \otimes_E (\Omega^\bullet_{\Nbf/E} \otimes_E V)^{\Nbf} \to \Omega^\bullet_{\Nbf/E} \otimes_E V \]
is an isomorphism in \(\Rep_{\alg}(\Nbf)\), in particular each term of \(\Omega^\bullet_{\Nbf/E} \otimes_E V\) is an injective object of \(\Rep_{\alg}(\Nbf)\).
It would follow from general formalism (see \cite[\href{https://stacks.math.columbia.edu/tag/013P}{Lemma 013P (2)}]{stacks-project}) that there exists a morphism \(R(\Nbf^{\bullet+1}) \to \Omega^\bullet_{\Nbf/E}\) which in degree \(0\) is the obvious isomorphism.
We will need an explicit formula however, in order to obtain a morphism of complexes
\[ C^\bullet_{\alg}(\Nbf, V) \longrightarrow \Omega^\bullet_{\Nbf/E} \otimes_E V \]
in \(\Rep_{\alg}(\Nbf \rtimes \Hbf)\), whenever \(\Hbf\) is a linear algebraic group acting on \(\Nbf\) and \(V\) is an object of \(\Rep_{\alg}(\Nbf \rtimes \Hbf)\).
Fortunately this formula is very simple:
\begin{align} \label{eq:Cpol_to_diff_forms}
  R(\Nbf^{n+1}) = R(\Nbf) \otimes_E \dots \otimes_E R(\Nbf) & \longrightarrow \Omega_{\Nbf/E}^n \\
  f_0 \otimes \dots \otimes f_n & \longmapsto f_0 \, df_1 \wedge \dots \wedge df_n. \nonumber
\end{align}
We may define \(\nfrak = \Lie \Nbf\) as the space of \(\Nbf\)-invariant (for left translation) derivations, i.e.\ elements of
\[ \Hom_{R(\Nbf)}(\Omega^1_{\Nbf/E}, R(\Nbf))^{\Nbf} \simeq \Hom_E((\Omega^1_{\Nbf/E})^{\Nbf}, E). \]
Recall that \(\Omega^1_{\Nbf/E}\) is generated as a \(R(\Nbf)\)-module by the elements \(df\) where \(f \in R(\Nbf)\).
As usual for \(X \in \Hom_{R(\Nbf)}(\Omega^1_{\Nbf/E}, R(\Nbf))\) we denote the image of \(df \in \Omega^1_{\Nbf/E}\) by \(X\) as \(df(X)\).
The Lie bracket
\[ [\cdot, \cdot]: \Hom_{R(\Nbf)}(\Omega^1_{\Nbf/E}, R(\Nbf))^2 \longrightarrow \Hom_{R(\Nbf)}(\Omega^1_{\Nbf/E}, R(\Nbf)) \]
is defined by the formula
\[ df([X,Y]) = d \left( df(Y) \right)(X) - d \left( df(X) \right)(Y) \]
and clearly maps \(\nfrak^2\) to \(\nfrak\).
We have isomorphisms\footnote{This uses the assumption \(\car E = 0\).} in \(\Rep_{\alg}(\Nbf)\)
\begin{align} \label{eq:Omega_as_forms_Lie}
  \Omega_{\Nbf/E}^n & \longrightarrow \Hom_E \left( \bigwedge^n \nfrak, R(\Nbf) \right) = \Hom_{R(\Nbf)} \left( \bigwedge^n \Hom_{R(\Nbf)}(\Omega^1_{\Nbf/E}, R(\Nbf)), R(\Nbf) \right) \\
  f_0 df_1 \wedge \dots \wedge df_n & \longmapsto \left( (X_1, \dots, X_n) \mapsto f_0 \sum_{\sigma \in S_n} \epsilon(\sigma) df_{\sigma(1)}(X_1) \dots df_{\sigma(n)}(X_n) \right). \nonumber
\end{align}
Via these isomorphisms the differential on \(\Omega^\bullet_{\Nbf/E}\) is given by the formula\footnote{This computation is well-known, at least in the setting of real Lie
groups.}
\begin{align*}
  \Hom_E \left( \bigwedge^n \nfrak, R(\Nbf) \right) & \longrightarrow \Hom_E \left( \bigwedge^{n+1} \nfrak, R(\Nbf) \right) \\
  \varphi & \longmapsto d\varphi
\end{align*}
\begin{align*}
  d\varphi(X_1, \dots, X_{n+1}) :={}&
  \sum_{i=1}^{n+1} (-1)^{i+1} d(\varphi(X_1, \dots, \widehat{X_i}, \dots, X_{n+1}))(X_i) \\
  &+ \sum_{1 \leq i < j \leq n+1} (-1)^{i+j} \varphi([X_i, X_j], X_1, \dots, \widehat{X_i}, \dots, \widehat{X_j}, \dots, X_{n+1}).
\end{align*}
For any object \(V\) of \(\Rep_{\alg}(\Nbf)\) the above isomorphism identifies \((\Omega^n_{\Nbf/E} \otimes_E V)^{\Nbf}\) with
\[ \Hom_E \left( \bigwedge^n \nfrak, (R(\Nbf) \otimes_E V)^{\Nbf} \right) \simeq \Hom_E \left( \bigwedge^n \nfrak, V \right) \]
by evaluation at \(1 \in \Nbf(E)\).
We recover the isomorphism from \((\Omega^\bullet_{\Nbf/E} \otimes_E V)^{\Nbf}\) to the Chevalley-Eilenberg complex \(C^\bullet(\nfrak, V)\) computing Lie algebra cohomology recalled in Section \ref{sec:Lie_alg_coh}\footnote{Strictly speaking we only defined the Chevalley-Eilenberg complex in the case where \(V\) is finite-dimensional, but of course it may be defined without this assumption.}.
More generally if \(\Hbf\) is a linear algebraic group over \(E\) acting on \(\Nbf\) and \(V\) is an object of \(\Rep_{\alg}(\Nbf \rtimes \Hbf)\) then the composition of \eqref{eq:Cpol_to_diff_forms} and \eqref{eq:Omega_as_forms_Lie} yields a morphism of complexes in \(\Rep_{\alg}(\Hbf)\)
\[ \Theta_{\Nbf}: C^\bullet_{\alg}(\Nbf, V) \to C^\bullet(\nfrak, V). \]

If \(\Nbf\) is connected and unipotent then the complex \(\Omega^\bullet_{\Nbf/E}\) is exact in positive degree because \(\Nbf\) is isomorphic, as a scheme over \(E\), to the affine space \(\A^{\dim \Nbf}_E\).
It follows that \(\Omega^\bullet_{\Nbf/E} \otimes_E V\) is a resolution of \(V\) by injective objects of \(\Rep_{\alg}(\Nbf)\) in this case, and so \(\Theta_{\Nbf}\) is a quasi-isomorphism.
By the usual arguments involving double complexes, if \(V^\bullet\) is a bounded below complex of objects of \(\Rep_{\alg}(\Nbf \rtimes \Hbf)\) then we have a quasi-isomorphism between bounded-below complexes of objects of \(\Rep_{\alg}(\Hbf)\)
\begin{equation} \label{eq:tot_pol_cocy_qis_CE}
  \Theta_{\Nbf}: \Tot^\bullet \left( C^\bullet_{\alg}(\Nbf, V^\bullet) \right) \longrightarrow \Tot^\bullet \left( \Hom_E \left( \bigwedge^\bullet \nfrak, V^\bullet \right) \right).
\end{equation}

The definitions and arguments in Section \ref{sec:expl_HS} have obvious analogues for the complexes \(C^\bullet_{\alg}(\Nbf, V)\), in particular the explicit version of the Hochschild-Serre spectral sequence (Theorem \ref{thm:Eilenberg_Zilber}) holds in this setting as well.
We conclude this section with the comparison of this result, the maps \(\Theta_{\Nbf}\) defined above and Lemma \ref{lem:HS_Lie_adhoc}.

\begin{lemm} \label{lem:Theta_normalized}
  For \(1 \leq j \leq n\) define
  \begin{align*}
    \tau_j: R(\Nbf^{n+1}) \otimes_E V & \longrightarrow R(\Nbf^{n+1}) \otimes_E V \\
    f_0 \otimes \dots \otimes f_n \otimes v & \longmapsto f_j \otimes f_0 \otimes \dots \widehat{f_j} \dots \otimes f_n \otimes v
  \end{align*}
  and define \(\tau_0 = \id\).
  On the subspace \(C^n_{\alg,\norm}(\Nbf, V)\) of normalized cochains, defined as in \eqref{eq:def_norm_cochains}, we have
  \[ \Theta_{\Nbf} = (-1)^j \Theta_{\Nbf} \circ \tau_j. \]
\end{lemm}
\begin{proof}
  Neither the tensor product with \(V\) nor the invariance condition under \(\Nbf\) play any role in the proof and so we omit both.
  Recall that \(\Theta_{\Nbf}\) is defined as the composition of \eqref{eq:Cpol_to_diff_forms}, denoted \(\alpha\) in this proof, and \eqref{eq:Omega_as_forms_Lie}, denoted \(\beta\) in this proof.
  (This composition is restricted to \(\Nbf\)-invariant forms, but the distinction is irrelevant here.)
  By definition the subspace of normalized cochains is included in the kernel of
  \begin{align*}
    \pi_i: R(\Nbf^{n+1}) & \longmapsto R(\Nbf^n) \\
    f_0 \otimes \dots \otimes f_n & \longmapsto f_i f_{i+1} \otimes f_0 \otimes \dots \otimes f_{i-1} \otimes f_{i+2} \otimes \dots \otimes f_n
  \end{align*}
  for any \(0 \leq i \leq n-1\).
  With the convention \(\tau_0 = \id\) we compute
  \[ (-1)^i d \circ \alpha \circ \pi_i = \alpha \circ \tau_i + \alpha \circ \tau_{i+1}.  \]
  It follows that on normalized cochains we have
  \[ \alpha = -\alpha \circ \tau_1 = \dots = (-1)^j \alpha \circ \tau_j. \]
  Applying \(\beta\) yields the desired formula.
\end{proof}

\begin{prop} \label{pro:comp_HS_alggp_Lie}
  Let \(\Nbf_1\) and \(\Nbf_2\) be unipotent algebraic groups over \(E\).
  Let \(\Hbf\) be a linear algebraic group over \(E\).
  Assume given an action of \(\Hbf\) on \(\Nbf_2\) and an action of \(\Nbf_2 \rtimes \Hbf\) on \(\Nbf_1\).
  For any object \(V\) of \(\Rep_{\fg,\alg}(\Nbf_1 \rtimes (\Nbf_2 \rtimes \Hbf))\) and any \(n \in \Z_{\geq 0}\) the following diagram commutes up to homotopy, where the bottom right arrow is the isomorphism of Lemma \ref{lem:HS_Lie_adhoc}.
  \[
    \begin{tikzcd}
      \bigoplus_{i+j=n} C^i_{\alg}(\Nbf_2, C^j_{\alg}(\Nbf_1, V)) \arrow[rr, "{\AW}"] \arrow[d, "{\Theta_{\Nbf_1}}"] && C^n_{\alg}(\Nbf_1 \rtimes \Nbf_2, V) \arrow[d, "{\Theta_{\Nbf_1 \rtimes \Nbf_2}}"] \\
      \bigoplus_{i+j=n} C^i_{\alg}(\Nbf_2, C^j(\nfrak_1, V)) \arrow[dr, "{\Theta_{\Nbf_2}}"] && C^n(\nfrak_1 \oplus \nfrak_2, V) \arrow[dl, "{\sim}"] \\
      & \bigoplus_{i+j=n} C^i(\nfrak_2, C^j(\nfrak_1, V)) &
    \end{tikzcd}
  \]
  In particular the same holds if we replace \(V\) by a bounded complex and take total complexes in the diagram.
\end{prop}
\begin{proof}
  It is enough to show that the restriction of both compositions to
  \[ C^i_{\pol,\norm}(\Nbf_2, C^j_{\pol,\norm}(\Nbf_1, V)) \]
  are equal.
  Recall from Theorem \ref{thm:Eilenberg_Zilber} that \(\AW\) maps this space to \(C^{i+j}_{\pol,\norm}(\Nbf_1 \rtimes \Nbf_2, V)\).
  For explicit computations it is convenient to denote
  \[ c = f_0 \otimes \dots \otimes f_i \otimes (g_0 \otimes \dots \otimes g_j \otimes v) \in R(\Nbf_2^{i+1}) \otimes_E (R(\Nbf_1^{j+1}) \otimes_E V). \]
  By definition we have
  \[ \AW(c) = (g_0 \otimes 1) \otimes \dots \otimes (g_{j-1} \otimes 1) \otimes (g_j \otimes f_0) \otimes (1 \otimes f_1) \otimes \dots \otimes (1 \otimes f_i) \otimes v \]
  in \(R((\Nbf_1 \rtimes \Nbf_2)^{n+1}) \otimes_E V\).
  Applying \(\Theta_{\Nbf_1 \rtimes \Nbf_2} \circ \tau_j\) and evaluating at \((X_1, \dots, X_n) \in (\nfrak_1 \oplus \nfrak_2)^n\) we obtain
  \[ (g_j \otimes f_0) \sum_{\sigma \in S_n} \epsilon(\sigma) d(g_0 \otimes 1)(X_{\sigma(1)}) \dots d(g_{j-1} \otimes 1)(X_{\sigma(j)}) d(1 \otimes f_1)(X_{\sigma(j+1)}) \dots d(1 \otimes f_i)(X_{\sigma(n)}) v. \]
  If we impose that \(X_k\) belongs to \(\nfrak_1\) for \(k \leq j'\) and to \(\nfrak_2\) for \(k>j'\), then the term corresponding to \(\sigma\) vanishes unless we have \(j'=j\) and \(\sigma \in S_j \times S_i\).
  More precisely we obtain
  \begin{align*}
    & \left( \Theta_{\Nbf_1 \rtimes \Nbf_2} \tau_j \AW(c) \right) (X_1, \dots, X_{j'}, Y_1, \dots, Y_{n-j'}) \\
    =&
    \begin{cases}
      \left( \Theta_{\Nbf_1} \tau_j \left( \Theta_{\Nbf_2}(c)(Y_1, \dots, Y_i) \right) \right) (X_1, \dots, X_j) & \text{ if } j'=j, \\
      0 & \text{ otherwise.}
    \end{cases}
  \end{align*}
  whenever \(X_1, \dots, X_{j'} \in \nfrak_1\) and \(Y_1, \dots, Y_{n-j'} \in \nfrak_2\).
  By Lemma \ref{lem:Theta_normalized} this implies
  \[ \Theta_{\Nbf_1 \rtimes \Nbf_2}(\AW(c)) = \Theta_{\Nbf_2} \Theta_{\Nbf_1} (c) \]
  for any \(c \in C^i_{\alg,\norm}(\Nbf_2, C^j_{\alg,\norm}(\Nbf_1, V))\)
\end{proof}

\section{A Nomizu-van Est theorem}
\label{sec:pf_vanEst}

\begin{theo} \label{thm:ladic_vanEst}
  Let \(\Nbf\) be a unipotent linear algebraic group over \(\Qell\).
  Let \(\Hbf\) be a linear algebraic group over \(\Qell\) acting on \(\Nbf\).
  Let \(K\) be a compact open subgroup of \((\Nbf \rtimes \Hbf)(\Qell)\).
  Denote \(N_K = K \cap \Nbf(\Qell)\) and \(H_K = K/N_K\), considered as a compact open subgroup of \(\Hbf(\Qell)\).
  Let \(\nfrak\) be the Lie algebra of \(\Nbf_E\).
  We have an isomorphism of composite functors
  \[ \begin{tikzcd}[column sep=7em, row sep=5em]
      D^b(\Rep_{\fg,\alg}((\Nbf \rtimes \Hbf)_E)) \arrow[r, "{R\Gamma_{\Lie}(\nfrak, -)}" above] \arrow[d, "{F}" left] & D^b(\Rep_{\fg,\alg}(\Hbf_E)) \arrow[d, "{F}" right] \\
      D^+(K, E) \arrow[r, "{R\Gamma(N_K, -)}" below] \arrow[ur, Rightarrow, shorten=16mm, "{\sim}" below right, "{\NvE_{\Nbf,\Hbf,K}}" above left] & D^+(H_K, E)
    \end{tikzcd} \]
  i.e.\ isomorphisms in \(D^+(H_K, E)\)
  \[ \NvE_{\Nbf, \Hbf, K}(V^\bullet): R\Gamma(N_K, F(V^\bullet)) \simeq F(R\Gamma_{\Lie}(\nfrak, V^\bullet)) \]
  for all objects \(V^\bullet\) of \(D^b(\Rep_{\fg, \alg}((\Nbf \rtimes \Hbf)_E))\), which are functorial in \(V^\bullet\).
  In the diagram \(F\) abusively means the composition of the obvious functor
  \[ D^b(\Rep_{\fg,\alg}((\Nbf \rtimes \Hbf)_E)) \to D^b(\Rep_{\fg,\cont}(K, E)) \]
  with the functor \(F\) of Proposition \ref{pro:F_to_Eke}, and similarly for \(\Hbf_E\) and \(H_K\).

  These isomorphisms enjoy the following compatibility properties.
  \begin{enumerate}
  \item The isomorphisms \(\NvE_{\Nbf,\Hbf,K}\) are compatible with restriction to a subgroup \(\Hbf'\) of \(\Hbf\).
  \item For \(K'\) any open subgroup of \(K\) and any object \(V^\bullet\) of \(D^b(\Rep_{\fg, \alg}((\Nbf \rtimes \Hbf)_E))\) the following diagram is commutative
    \[
      \begin{tikzcd}
        \res_{H_{K'}}(R\Gamma(N_K, F(V^\bullet))) \arrow[dr, "{\sim}" above right, "{\res_{H_{K'}}(\NvE_{\Nbf,\Hbf,K}(V^\bullet))}" below left] \arrow[rr, "{r_{K,N_K,K'}(F(V^\bullet))}"] && R\Gamma(N_{K'}, F(V^\bullet)) \arrow[dl, "{\sim}" above left, "{\NvE_{\Nbf,\Hbf,K'}(V^\bullet)}" below right] \\
        & F(R\Gamma_{\Lie}(\nfrak, V^\bullet)) &
      \end{tikzcd}
    \]
  \item Assume that we have a decomposition \(\Nbf = \Nbf_1 \rtimes \Nbf_2\) which is preserved by the action of \(\Hbf\), in other words we have \(\Nbf \rtimes \Hbf = \Nbf_1 \rtimes (\Nbf_2 \rtimes \Hbf)\).
    Denote \(N_{1,K} = K \cap \Nbf_1(\Qell)\), \(N_{2,K} = (K \cap (\Nbf_1 \rtimes \Nbf_2)(\Qell))/N_{1,K}\) (understood as a compact open subgroup of \(\Nbf_2(\Qell)\)) and \(U_K = K/N_{1,K}\) (understood as a compact open subgroup of \((\Nbf_2 \rtimes \Hbf)(\Qell)\)).
    Then for any object \(V^\bullet\) of \(D^b(\Rep_{\fg, \alg}((\Nbf \rtimes \Hbf)_E))\) the following diagram is commutative
    \[
      \begin{tikzcd}[column sep=6em]
        R\Gamma(N_K, F(V^\bullet)) \arrow[r, "{\sim}" below, "{\NvE_{\Nbf,\Hbf,K}}" above] \arrow[d, "{\sim}"] & F(R\Gamma_{\Lie}(\nfrak, V^\bullet)) \arrow[dd, "{\sim}" right] \\
        R\Gamma(N_{2,K}, R\Gamma(N_{1,K}, F(V^\bullet))) \arrow[d, "{\sim}" right, "{R\Gamma(N_{2,K}, \NvE_{\Nbf_1, \Nbf_2 \rtimes \Hbf, K}(V^\bullet))}" left] & \\
        R\Gamma(N_{2,K}, F(R\Gamma_{\Lie}(\nfrak_1, V^\bullet))) \arrow[r, "{\sim}" below, "{\NvE_{\Nbf_2, \Hbf, U_K}}" above] & F(R\Gamma_{\Lie}(\nfrak_2, R\Gamma_{\Lie}(\nfrak_1, V^\bullet)))
      \end{tikzcd}
    \]
    Here the left vertical isomorphism is given by Lemma \ref{lem:HS_Lie_adhoc} and the top left vertical isomorphism is the usual composition of derived functors.
  \end{enumerate}
\end{theo}
\begin{proof}
  Endow \(\Nbf\) with a filtration \((\Nbf_i)_{0 \leq i \leq m+1}\) (e.g.\ the upper central series) and weights \(\ul{w} = (w_A, w_m, \dots, w_0)\) making it a weighted filtered unipotent algebraic group in the sense of Section \ref{sec:char_0_coeff}.
  Denote \(N_i = N \cap \Nbf_i(\Qell)\).
  Define a filtration \(\Fil^\bullet V^\bullet\) on \(V^\bullet\) by \(\Fil^j V^\bullet = 0\) for \(j<0\) and \(\Fil^j V^n = (V^n/\Fil^{j-1} V^n)^{\Nbf_E}\) for \(j \geq 0\).
  By Proposition \ref{pro:cont_coh_as_pol_ell_inv} we have an isomorphism in \(D^+(H_K, E)\)
  \begin{equation} \label{eq:step1_pf_vanEst}
    R\Gamma(N_K, F(V^\bullet)) \simeq F \left( \Tot^\bullet(\Fil^d_{\ul{w}} C^\bullet_{\pol}(\Nbf(\Qell), V^\bullet)) \right)
  \end{equation}
  for any large enough integer \(d\).
  For large enough \(d\) we also have a quasi-isomorphism
  \[ \Tot^\bullet \left( \Fil^d_{\ul{w}} C^\bullet_{\pol}(\Nbf(\Qell), V^\bullet) \right) \subset \Tot^\bullet \left( C^\bullet_{\pol}(\Nbf(\Qell), V^\bullet) \right) \]
  thanks to Lemma \ref{lem:alg_N_Fil_Cpol_stabilizes}.
  We also have a quasi-isomorphism \eqref{eq:tot_pol_cocy_qis_CE}
  \[ \Tot^\bullet \left( C^\bullet_{\pol}(\Nbf(\Qell), V^\bullet) \right) \xrightarrow{\sim} R\Gamma_{\Lie}(\nfrak, V^\bullet). \]
  Composing the two gives an isomorphism in \(D^b(\Rep_{\fg,\cont}(H_K, E))\)
  \[ \Tot^\bullet \left( \Fil^d_{\ul{w}} C^\bullet_{\pol}(\Nbf(\Qell), V^\bullet) \right) \to R\Gamma_{\Lie}(\nfrak, V^\bullet). \]
  Applying the functor \(F\) and composing with \eqref{eq:step1_pf_vanEst} gives \(\NvE_{\Nbf, \Hbf, K}(V^\bullet)\).
  It is easy to check that this morphism does not depend on the choices made (filtration, weights, \(d\), model over \(\Ocal_E\) \dots), by cofinality arguments (e.g.\ \(\varinjlim_d \Fil^d_{\ul{w}} C^n_{\pol}(\Nbf(\Qell), V^m) = C^n_{\pol}(\Nbf(\Qell), V^m)\)).
  Using similar arguments one can check that \(\NvE_{\Nbf,\Hbf,K}(V^\bullet)\) is functorial in \(V^\bullet\).
  The details of these relatively formal arguments are omitted here.

  The proof of the first compatibility is very formal.
  
  The second compatibility property, where \(K'\) is an open subgroup of \(K\), follows from Lemma \ref{lem:res_iso_after_inv_ell}, also using arguments as in the proof of this lemma.

  The second compatibility property follows from Proposition \ref{pro:comp_HS_alggp_Lie}.
\end{proof}

\begin{rema}
  In applications it is often the case that \(\Hbf\) is reductive, so that the abelian category \(\Rep_{\fg,\alg}(\Hbf_E)\) is semi-simple and so \(D^b(\Rep_{\fg,\alg}(\Hbf_E))\) is abelian and its objects decompose as direct sums of complexes concentrated in one degree.
  In fact Theorem \ref{thm:ladic_vanEst} is meant to be applied in the case where \(\Nbf \rtimes \Hbf\) is a parabolic subgroup of a connected reductive group \(\Gbf\) over \(\Qell\) and we pre-compose with the restriction functor
  \[ D^b(\Rep_{\fg,\alg}(\Gbf_E)) \to D^b(\Rep_{\fg,\alg}((\Nbf \rtimes \Hbf)_E)). \]
  In this case we may start from an irreducible representation of \(\Gbf_E\) and Lie algebra cohomology was explicitly computed by Kostant \cite[Theorem 5.14]{Kostant_Liealgcoh}.
\end{rema}

For the application that motivates this work we need the following slightly more complicated consequence of Theorem \ref{thm:ladic_vanEst}.

\begin{coro} \label{cor:ladic_vanEst_overQ}
  Let \(\Nbf\) be a unipotent linear algebraic group over \(\Q\).
  Let \(\Hbf\) be a linear algebraic group over \(\Q\) acting on \(\Nbf\).
  Let \(K\) be a compact open subgroup of \(\Nbf(\A_f)\).
  Denote \(N_K = K \cap \Nbf(\A_f)\) and \(H_K = K/N_K\), considered as a compact open subgroup of \(\Hbf(\A_f)\).
  Let \(\nfrak\) be the Lie algebra of \(\Nbf_E\).
  We consider algebraic representations of \(\Hbf_E\) as continuous representations of \(K\) over \(E\) via the projection \(\Hbf(\A_f) \to \Hbf(\Qell)\).
  We have an isomorphism of composite functors
  \[ \begin{tikzcd}[column sep=7em, row sep=5em]
      D^b(\Rep_{\fg,\alg}((\Nbf \rtimes \Hbf)_E)) \arrow[r, "{R\Gamma_{\Lie}(\nfrak, -)}" above] \arrow[d, "{F}" left] & D^b(\Rep_{\fg,\alg}(\Hbf_E)) \arrow[d, "{F}" right] \\
      D^+(K, E) \arrow[r, "{R\Gamma(N_K, -)}" below] \arrow[ur, Rightarrow, shorten=16mm, "{\sim}" below right, "{\NvE_{\Nbf,\Hbf,K}}" above left] & D^+(H_K, E).
    \end{tikzcd} \]
  It satisfies compatibilities similar to the ones in Theorem \ref{thm:ladic_vanEst} (details left to the reader).
\end{coro}
\begin{proof}
  Denote by \(\A_f^{(\ell)}\) the restricted product of \(\Qp\) for primes \(p \neq \ell\), so that we have \(\A_f = \Qell \times \A_f^{(\ell)}\).
  Denote \(N_{K,\ell} = N_K \cap \Nbf(\Qell)\) and \(N_K^{(\ell)} = N_K \cap \Nbf(\A_f^{(\ell)})\).
  Denote \(K^{(\ell)} = K \cap (\Nbf \rtimes \Hbf)(\A_f^{(\ell)}\) and \(K_\ell = K/K^{(\ell)}\), considered as a compact open subgroup of \((\Nbf \rtimes \Hbf)(\Qell)\).
  We have \(N_K = N_{K,\ell} \times N_K^{(\ell)}\) and \(N_K^{(\ell)}\) is an inverse limit of finite groups having cardinality coprime to \(\ell\).
  It follows that the functor of \(N_K^{(\ell)}\)-invariants, from \((S_K^{\N}, (\Ocal_E)_\bullet)\)-modules to \((S_{K/N_K^{(\ell)}}^{\N}, (\Ocal_E)_\bullet)\)-modules, is exact.
  We get an isomorphism of functors
  \[
    \begin{tikzcd}[row sep=4em]
      D^+(K_\ell, E) \arrow[rd, "{\inf}" below left] \arrow[rr, "{\inf}" above] & \arrow[d, Rightarrow, shorten=6mm]& D^+(K/N_K^{(\ell)}, E) \\
      & D^+(K, E) \arrow[ru, "{R\Gamma(N_K^{(\ell)}, -)}" below right] &
    \end{tikzcd}
  \]
  where each map labelled \(\inf\) is an inflation map.
  We also have a morphism of functors
  \[
    \begin{tikzcd}[column sep=4em, row sep=4em]
      D^+(K_\ell, E) \arrow[r, "{R\Gamma(N_{K,\ell}, -)}"] \arrow[d, "{\inf}"] & D^+(K_\ell/N_{K,\ell}, E) \arrow[d, "{\inf}"] \arrow[Rightarrow, ld, shorten=8mm] \\
      D^+(K/N_K^{(\ell)}, E) \arrow[r, "{R\Gamma(N_{K,\ell}, -)}"] & D^+(H_K, E)
    \end{tikzcd}
  \]
  which is easily seen to be an isomorphism using the isomorphisms of functors \(r_{K_\ell, N_{K,\ell}, N_{K,\ell}}\) and \(r_{K/N_K^{(\ell)}, N_{K,\ell}, N_{K,\ell}}\) defined in \eqref{eq:res_RGamma}.
  Combining these two isomorphisms of functors with the isomorphism of functors
  \[ R\Gamma(N_K, -) \simeq R\Gamma(N_{K,\ell}, R\Gamma(N_K^{(\ell)}, -)) \]
  reduces the corollary to Theorem \ref{thm:ladic_vanEst}.
\end{proof}

\begin{rema} \label{rem:dont_understand_Pink}
  As mentioned in the introduction, the motivation for proving Corollary \ref{cor:ladic_vanEst_overQ} is to prove an analogue in triangulated categories, rather than just for cohomology groups, of \cite[Theorem 5.3.1]{Pink_ladic_Shim}.
  One might wonder why the proof of Corollary \ref{cor:ladic_vanEst_overQ} does not use, as in \S 5.2 loc.\ cit., an argument by restriction from \(K \cap \Nbf(\A_f)\) (denoted \(K_W\) loc.\ cit.) to \(K \cap \Nbf(\Q)\) (denoted \(\Gamma_W\) loc.\ cit.).
  Using the notation loc.\ cit., I do not understand how this argument by restriction can possibly yield a \(H_Q/H_C\)-equivariant isomorphism, because there seems to be no subgroup of \(H_Q\) surjecting onto \(H_Q/H_C\) and whose intersection with \(K_W\) is \(\Gamma_W\).
  Using our notation, the issue is roughly that the action of \(K \cap \Hbf(\A_f)\) on \(K \cap \Nbf(\A_f)\) does not preserve \(K \cap \Nbf(\Q)\) in general.
\end{rema}

\newpage

\printbibliography

@article{Lazard_gpanp,
 Author = {Lazard, Michel},
 Title = {Groupes analytiques {\(p\)}-adiques},
 FJournal = {Publications Math{\'e}matiques},
 Journal = {Publ. Math., Inst. Hautes {\'E}tud. Sci.},
 Volume = {26},
 Pages = {389--603},
 Year = {1965}
}

@unpublished{Zhu_IHorth,
author = {Zhu, Yihang},
title = {The stabilization of the Frobenius--Hecke traces on the intersection cohomology of orthogonal Shimura varieties},
note = {\url{https://arxiv.org/abs/1801.09404}}
}

@article{GoreskyHarderMacPherson_weighted,
    author = {Goresky, M. and Harder, G. and MacPherson, R.},
     title = {Weighted cohomology},
   journal = {Invent. Math.},
  fjournal = {Inventiones Mathematicae},
    volume = {116},
      year = {1994},
    number = {1-3},
     pages = {139--213}
}

@article{vanEst,
 Author = {van Est, W. T.},
 Title = {A generalization of the {Cartan}-{Leray} spectral sequence. {I}, {II}},
 FJournal = {Nederlandse Akademie van Wetenschappen. Proceedings. Series A. Indagationes Mathematicae},
 Journal = {Nederl. Akad. Wet., Proc., Ser. A},
 Volume = {61},
 Pages = {399--405, 406--413},
 Year = {1958},
}

@article{Nomizu,
 Author = {Nomizu, Katsumi},
 Title = {On the cohomology of compact homogeneous spaces of nilpotent {Lie} groups},
 FJournal = {Annals of Mathematics. Second Series},
 Journal = {Ann. Math. (2)},
 Volume = {59},
 Pages = {531--538},
 Year = {1954},
 DOI = {10.2307/1969716},
}

@misc{stacks-project,
    shorthand    = {Stacks},
    author       = {The {Stacks Project Authors}},
    title        = {\textit{Stacks Project}},
    howpublished = {\url{https://stacks.math.columbia.edu}},
    year         = {2022}
}

@misc{DemazureGabriel,
 Author = {Demazure, Michel and Gabriel, Pierre},
 Title = {Groupes alg{\'e}briques. {Tome} {I}: {G{\'e}om{\'e}trie} alg{\'e}brique. {G{\'e}n{\'e}ralit{\'e}s}. {Groupes} commutatifs. {Avec} un appendice `{Corps} de classes local' par {Michiel} {Hazewinkel}},
 Year = {1970},
 Language = {English},
 HowPublished = {Paris: {Masson} et {Cie}, {\'E}diteur; {Amsterdam}: {North}-{Holland} {Publishing} {Company}. xxvi, 700 p. (1970).}
}

@book{KargapolovMerzljakov_GTM,
 Author = {Kargapolov, M. I. and Merzljakov, Ju. I.},
 Title = {Fundamentals of the theory of groups. {Transl}. from the 2nd {Russian} ed. by {Robert} {G}. {Burns}},
 FSeries = {Graduate Texts in Mathematics},
 Series = {Grad. Texts Math.},
 Volume = {62},
 Year = {1979},
 Publisher = {Springer, Cham},
 Language = {English}
}

@article{Malcev_nilp_tors_free,
 Author = {Mal'tsev, A. I.},
 Title = {Torsion-free nilpotent groups},
 FJournal = {Izvestiya Akademii Nauk SSSR. Seriya Matematicheskaya},
 Journal = {Izv. Akad. Nauk SSSR, Ser. Mat.},
 Volume = {13},
 Pages = {201--212},
 Year = {1949},
 Language = {Russian}
}

@Book{Verdier_catder,
 Author = {Jean-Louis {Verdier}},
 Title = {{Des cat\'egories d\'eriv\'ees des cat\'egories ab\'eliennes}},
 FJournal = {{Ast\'erisque}},
 Journal = {{Ast\'erisque}},
 Volume = {239},
 Pages = {ix + 253},
 Year = {1996},
 Publisher = {Paris: Soci\'et\'e Math\'ematique de France},
 Language = {French},
}

@Article{Hochschild_cohalglingps,
 Author = {G. {Hochschild}},
 Title = {{Cohomology of algebraic linear groups}},
 FJournal = {{Illinois Journal of Mathematics}},
 Journal = {{Ill. J. Math.}},
 Volume = {5},
 Pages = {492--519},
 Year = {1961},
 Publisher = {Duke University Press, Durham, NC},
 Language = {English},
}

@phdthesis{Rubio_these,
  TITLE = {{Homologie effective des espaces de lacets it{\'e}r{\'e}s : un
logiciel}},
  AUTHOR = {Rubio-Garcia, Julio},
  URL = {https://tel.archives-ouvertes.fr/tel-00339304},
  SCHOOL = {{Universit{\'e} Joseph-Fourier - Grenoble I}},
  YEAR = {1991},
  MONTH = Oct,
  TYPE = {Theses},
  HAL_ID = {tel-00339304},
  HAL_VERSION = {v1},
}

@Article{EilenbergMacLane2,
 Author = {Samuel {Eilenberg} and Saunders {MacLane}},
 Title = {{On the groups \(H(\Pi,n)\). II}},
 FJournal = {{Annals of Mathematics. Second Series}},
 Journal = {{Ann. Math. (2)}},
 Volume = {60},
 Pages = {49--139},
 Year = {1954},
 Language = {English},
 DOI = {10.2307/1969702},
}

@Article{EilenbergMacLane1,
 Author = {Samuel {Eilenberg} and Saunders {MacLane}},
 Title = {{On the groups \(H(\Pi,n)\). I}},
 FJournal = {{Annals of Mathematics. Second Series}},
 Journal = {{Ann. Math. (2)}},
 Volume = {58},
 Pages = {55--106},
 Year = {1953},
 Language = {English},
 DOI = {10.2307/1969820},
}

@Book{KashiwaraShapira_catandsheaves,
 Author = {Masaki {Kashiwara} and Pierre {Schapira}},
 Title = {{Categories and sheaves}},
 FJournal = {{Grundlehren der Mathematischen Wissenschaften}},
 Journal = {{Grundlehren Math. Wiss.}},
 Volume = {332},
 Pages = {x + 497},
 Year = {2006},
 Publisher = {Berlin: Springer},
 Language = {English},
 DOI = {10.1007/3-540-27950-4},
}

@book{Serre_corpsloc,
    AUTHOR = {Serre, Jean-Pierre},
     TITLE = {Corps locaux},
    SERIES = {Publications de l'Universit\'{e} de Nancago, No. VIII},
      NOTE = {Deuxi\`eme \'{e}dition},
 PUBLISHER = {Hermann, Paris},
      YEAR = {1968},
     PAGES = {245},
   MRCLASS = {12BXX (14GXX)},
  MRNUMBER = {0354618},
}

@book{MorelBook,
    author = {Morel, Sophie},
     title = {On the cohomology of certain noncompact {S}himura varieties},
    series = {Annals of Mathematics Studies},
    volume = {173},
      note = {With an appendix by Robert Kottwitz},
 publisher = {Princeton University Press, Princeton, NJ},
      year = {2010},
     pages = {xii+217}
}

@article{MorelSiegel1,
    author = {Morel, Sophie},
     title = {Complexes pond\'er\'es sur les compactifications de
              {B}aily-{B}orel: le cas des vari\'et\'es de {S}iegel},
   journal = {J. Amer. Math. Soc.},
    volume = {21},
      year = {2008},
    number = {1},
     pages = {23--61}
}

@incollection {Ekedahl_adic,
    author = {Ekedahl, Torsten},
     title = {On the adic formalism},
 booktitle = {The {G}rothendieck {F}estschrift, {V}ol.\ {II}},
    series = {Progr. Math.},
    volume = {87},
     pages = {197--218},
 publisher = {Birkhäuser Boston, Boston, MA},
      year = {1990}
}

@book {SGA4-1,
     title = {Théorie des topos et cohomologie étale des schémas. {T}ome 1:
              {T}héorie des topos},
    series = {Lecture Notes in Mathematics, Vol. 269},
      note = {Séminaire de Géométrie Algébrique du Bois-Marie 1963--1964
              (SGA 4),
              Dirigé par M. Artin, A. Grothendieck, et J. L. Verdier. Avec
              la collaboration de N. Bourbaki, P. Deligne et B. Saint-Donat},
 publisher = {Springer-Verlag, Berlin-New York},
      year = {1972},
     pages = {xix+525}
}

@article {Pink_ladic_Shim,
    author = {Pink, Richard},
     title = {On {$l$}-adic sheaves on {S}himura varieties and their higher
              direct images in the {B}aily-{B}orel compactification},
   journal = {Math. Ann.},
    volume = {292},
      year = {1992},
    number = {2},
     pages = {197--240},
       doi = {10.1007/BF01444618},
}

@article {Kostant_Liealgcoh,
    author = {Kostant, Bertram},
     title = {Lie algebra cohomology and the generalized {B}orel-{W}eil
              theorem},
   journal = {Ann. of Math. (2)},
    volume = {74},
      year = {1961},
     pages = {329--387},
       doi = {10.2307/1970237},
}

\end{document}